\documentclass[12pt,leqno]{article}
\usepackage{amsmath,amssymb,amscd,latexsym}
\usepackage[all]{xy}
\newcommand{\C}{{\mathbb{C}}}
\newcommand{\F}{{\mathbb{F}}}
\newcommand{\Pa}{{\mathbb{P}}}
\newcommand{\Q}{{\mathbb{Q}}}
\newcommand{\R}{{\mathbb{R}}}
\newcommand{\Z}{{\mathbb{Z}}}
\newcommand{\complex}{\mathrm{complex}}
\newcommand{\ddet}{\mathrm{det}}
\renewcommand{\div}{\mathrm{div}\,}
\newcommand{\id}{\mathrm{id}}
\newcommand{\ord}{\mathrm{ord}}
\newcommand{\real}{\mathrm{real}}
\newcommand{\spec}{\mathrm{spec}\,}
\newcommand{\vol}{\mathrm{vol}\,}
\newcommand{\Fr}{\mathrm{Fr}}
\newcommand{\RRe}{\mathrm{Re}\,}
\newcommand{\Tr}{\mathrm{Tr}}
\newcommand{\Ch}{{\mathcal C}}
\newcommand{\Dh}{{\mathcal D}}
\newcommand{\tD}{\tilde{\Dh}}
\newcommand{\Fh}{{\mathcal F}}
\newcommand{\Nh}{{\mathcal N}}
\newcommand{\Oh}{{\mathcal O}}

\newcommand{\tU}{\tilde{U}}
\newcommand{\Wh}{\mathcal{W}}
\newcommand{\ea}{\mathfrak{a}}
\newcommand{\ed}{\mathfrak{d}}
\newcommand{\eo}{\mathfrak{o}}
\newcommand{\ep}{\mathfrak{p}}
\newcommand{\ohne}{\smallsetminus}
\newcommand{\betr}[1]{|#1|}
\newcommand{\Betr}[1]{|\!|#1|\!|}
\newcommand{\silo}{\stackrel{\sim}{\longrightarrow}}
\newcommand{\tei}{\, | \,}
\newcommand{\hullet}{\raisebox{0.05cm}{$\scriptscriptstyle \bullet$}}
\newcommand{\halb}{\frac{1}{2}}

\newcommand{\rprod}[2]{\mbox{$\textstyle\prod \hspace{-4mm} {\rule[-3pt]{2.5mm}{0.2mm}}$}\,^{#1}_{#2}}
\newtheorem{theorem}{Theorem}[section]
\newtheorem{lemma}[theorem]{Lemma}
\newtheorem{utheorem}{Theorem}[theorem]
\newtheorem{ulemma}[utheorem]{Lemma}
\newtheorem{prop}[theorem]{Proposition}
\newtheorem{uprop}[utheorem]{Proposition}
\newtheorem{defn}[theorem]{Definition}
\newtheorem{cor}[theorem]{Corollary}
\newtheorem{ucor}[utheorem]{Corollary}

\newtheorem{remark}[theorem]{Remark}
\newenvironment{rem}{\bigskip\noindent{\bf Remark}}{\mbox{}\bigskip}
\newenvironment{rems}{\bigskip\noindent{\bf Remarks}}{\mbox{}\bigskip}
\newenvironment{proofof}{\bigskip\noindent{\bf Proof of}}{\mbox{}\hspace*{\fill}$\Box$ \bigskip}
\newenvironment{proof}{\bigskip\noindent{\bf Proof}}{\mbox{}\hspace*{\fill}$\Box$ \bigskip}
\newcommand{\beweisende}{\hspace*{\fill}$\Box$}
\begin{document}
%%%%%%%%%%%%%%%%%%%%%%%%%%%%%%%%%%%%%%%%%%%%%%%%%%%%%%%%%%
\title{Two-variable zeta functions and regularized products\\[0.3cm]
%{\normalsize \it Dedicated to Kazuya Kato}
} 
%\author{Christopher \surname{Deninger} \email{deninge@math.uni-muenster.de}}
\author{Christopher Deninger}
%\institute{Mathematisches Institut, WWU M\"unster, Einsteinstr. 62, 48149 M\"unster, Germany}
%\author{Wilhelm \surname{Singhof} \email{singhof@cs.uni-duesseldorf.de}}
%\institute{Mathematisches Institut, Universit\"at D\"usseldorf, Universit\"atsstr. 1, 40225 D\"usseldorf, Germany}
%\runningauthor{Christopher Deninger, Wilhelm Singhof}
%\runningtitle{Real polarizable Hodge structures arising from foliations}
\date{\ }
\maketitle
%%%%%%%%%%%%%%%%%%%%%%%%%%%%%%%%%%%%%%%%%%%%%%%%%%%%%%%%%%%
%\setcounter{section}{1}
%\input{sec1}
\section{Introduction}
\label{sec:1}

In his paper \cite{P} Pellikaan studied an interesting two-variable zeta function for algebraic curves over finite fields. Using notions from Arakelov theory of arithmetic curves, van der Geer and Schoof were led to introduce an analogous zeta function for number fields \cite{GS}.

In \cite{LR} Lagarias and Rains investigated this two-variable zeta function thoroughly for the special case of the rational number field. They also made some comments on the general case.

In earlier work we introduced a conjectural cohomological formalism to express Dedekind and more general zeta functions as regularized determinants of a certain operator $\Theta$ on cohomology. In this framework it is not unreasonable to assume that the second variable $w$ of the two-variable zeta function corresponds to an operator $\Theta_w$ depending on $w$. These heuristics which are explained in the last section suggest a formula for the two-variable zeta function as a regularized product.

The main contribution of the paper is to prove this formula for the two-variable zeta function of any number field, Theorem \ref{t52}. Our method is based on a powerful criterion of Illies for zeta-regularizability \cite{I1}, \cite{I2}. We refer to section \ref{sec:5} for a short review of the relevant facts from the theory of regularization.

We also treat the much easier case of curves over finite fields. For number fields, our approach requires us to determine the asymptotic behaviour for $\RRe s \to \infty$ of certain oscillatory integrals over spaces of lattices $\Gamma$. The function to be integrated is $a^{-s}_{\Gamma}$ where $a_{\Gamma}$ is the minimal length among the non-zero vectors in $\Gamma$. This is an interesting problem already for real quadratic fields in which case Don Zagier found a solution. The general case is treated in section \ref{sec:4}. 

The treatment in \cite{GS} and \cite{LR} of the two-variable zeta function for general number fields is somewhat brief. Also, the precise analogy with Pellikaan's original zeta function is not written down. In the first two sections we therefore give a more detailed exposition of these topics.

I would like to thank Don Zagier very much for his help in the real quadratic case which was a great inspiration for me. I am also grateful to Eva Bayer and Georg Illies for useful remarks and to the CRM in Montreal for its support.

\section{Background on two-variable zeta functions for curves over finite fields}
\label{sec:2}

Consider an algebraic curve $X$ over the finite field $\F_q$ with $q = p^r$ elements. Let $\betr{X}$ be the set of closed points of $X$ and for $x \in \betr{X}$ set \\
$\deg x = (\F_q (x) : \F_q)$. The zeta function of $X$ is defined by the Euler product
\[
Z_X (T) = \prod_{x \in \betr{X}} (1 - T^{\deg x})^{-1} \quad \mbox{in} \; \Z [\betr{T}] \; .
\]
For a divisor $D = \sum_{x \in \betr{X}} n_x \cdot x$ with $n_x \in \Z$ we set $\deg D = \sum n_x \deg x$. Then we have
\begin{equation}
  \label{eq:1}
  Z_X (T) = \sum_{D \ge 0} T^{\deg D}
\end{equation}
where the sum runs over all effective divisors i.e. those with $n_x \ge 0$ for all $x \in \betr{X}$. Let $CH^1 (X)$ denote the divisor class group of $X$ and for $\Dh = [D]$ set
\[
h^i (\Dh) := h^i (D) = \dim H^i (X , \Oh (D)) \; .
\]
Summing over divisor classes in (\ref{eq:1}), one gets the formula:
\begin{equation}
  \label{eq:2}
  Z_X (T) = \sum_{\Dh} \frac{q^{h^0 (\Dh)}-1}{q-1} T^{\deg \Dh} \; .
\end{equation}
Here it is enough to sum over $\Dh$'s with $\deg \Dh := \deg D \ge 0$. In \cite{P} \S\,3 Pellikaan had the idea to replace $q$ by a variable $u$ in this formula. His two-variable zeta function is defined by
\begin{equation}
  \label{eq:3}
  Z_X (T,u) = \sum_{\Dh} \frac{u^{h^0 (\Dh)}-1}{u-1} T^{\deg \Dh}  \; .
\end{equation}
Reconsidering classical proofs he obtained the following properties in the case where $X$ is smooth projective and geometrically irreducible:
\begin{gather}
  \label{eq:4}
  Z_X (T,u) = \frac{P_X (T,u)}{(1 -T) (1-uT)} \quad \mbox{with} \; P_X (T,u) \in \Z [T,u] \\
P_X (T,u) = \sum^{2g}_{i=0} P_i (u) T^i \quad \mbox{with} \; P_i (u) \in \Z [u] \; , \; \mbox{where} \label{eq:5} \\
  \label{eq:6}
  P_0 (u) = 1 , P_{2g} (u) = u^g \; , \; \deg P_i (u) \le 1 + \frac{i}{2} \quad \mbox{and $g$ is the genus of} \; X \; .
\end{gather}
The two-variable zeta function enjoys the functional equation
\begin{equation}
  \label{eq:7}
  Z_X (T,u) = u^{g-1} T^{2g-2} Z_X \left( \frac{1}{Tu} , u \right) \; .
\end{equation}
In terms of the $P_i (u)$ it reads:
\begin{equation}
  \label{eq:8}
  P_{2g-i} (u) = u^{g-i} P_i (u) \; .
\end{equation}
For example, for $X = \Pa^1$ one has $P_X (T,u) = 1$ and for $X$ an elliptic curve $P_X (T,u) = 1 + (\betr{X (\F_q)} - 1 - u) T + uT^2$. 

Recently Naumann \cite{N} proved the interesting fact that the polynomial $P_X (T,u)$ is irreducible in $\C [T,u]$. 

In \cite{GS} \S\,7, van der Geer and Schoof consider the following variant of Pellikaan's zeta function. They show that for complex $s$ and $t$ in $\RRe s < 0 , \RRe t < 0$ the series 
\begin{equation}
  \label{eq:9}
  \zeta^{GS}_X (s,t) = \sum_{\Dh \in CH^1 (X)} q^{sh^0 (\Dh) + t h^1 (\Dh)}
\end{equation}
defines a holomorphic function with a meromorphic continuation to $\C \times \C$. The explicit relation with $Z_X (T,u)$ is not stated in \cite{GS}, so we give it here:

\begin{prop}
  \label{t21}
  \begin{eqnarray*}
    \zeta^{GS}_X (s,t) & = & (q^{s+t} -1) q^{t (g-1)} Z_X (q^{-t} , q^{s+t}) \\
& = & (q^{s+t} -1) q^{s (g-1)} Z_X (q^{-s} , q^{s+t}) \; .
  \end{eqnarray*}
\end{prop}

\begin{proof}
  Using the Riemann--Roch theorem one obtains, c.f. \cite{GS} proof of prop.~5:
\[
\zeta^{GS}_X (s,t) = q^{t (g-1)} \sum_{0 \le \deg \Dh \le 2 g - 2} q^{(s + t) h^0 (\Dh)} q^{-t \deg \Dh} + h \left( \frac{q^{sg}}{1-q^s} + \frac{q^{tg}}{1 - q^t} \right) \; .
\]
Here $h$ is the order of $CH^1 (X)^0$, the group of degree zero divisor classes on $X$. This gives the meromorphic continuation to $\C \times \C$. On the other hand according to \cite{P}, p. 181 setting $u = q^{s+t} , T = q^{-t}$ we have:
\[
(q^{s+t} -1) Z_X (q^{-t} , q^{s+t}) = \sum_{0 \le \deg \Dh \le 2g - 2} q^{(s+t) h^0 (\Dh)} q^{-t \deg \Dh} + h \left( \frac{q^{sg + t (1-g)}}{1 - q^s} - \frac{1}{1 - q^{-t}} \right) \; .
\]
This implies the first equality in the proposition. The second follows from the functional equation (\ref{eq:7}) of $Z_X (T,u)$. 
\end{proof}
\bigskip

In particular the second relation in the proposition shows that for $s + t = 1$ we have
\begin{equation}
  \label{eq:10}
  \zeta^{GS}_X (s , 1-s) = (q-1) q^{s (g-1)} \zeta_X (s) \quad \mbox{where} \; \zeta_X (s) = Z_X (q^{-s})
\end{equation}
as stated in \cite{GS} proposition 5. Note that for $\zeta^{GS}_X (s,t)$ the functional equation takes the simple form:
\begin{equation}
  \label{eq:11}
  \zeta^{GS}_X (s,t) = \zeta^{GS}_X (t,s) \; .
\end{equation}
In the number field case, Lagarias and Rains introduced the substitution $t = w-s$. Thus we define here as well
\begin{equation}
  \label{eq:12}
  \zeta_X (s,w) = \zeta^{GS}_X (s,w-s) = (q^w -1) q^{-s (1-g)} Z_X (q^{-s} , q^w) \; .
\end{equation}
This meromorphic function of $s$ and $w$ satisfies the functional equation
\begin{equation}
  \label{eq:13}
  \zeta_X (s,w) = \zeta_X (w-s ,w) 
\end{equation}
and for $w = 1$ we have:
\begin{equation}
  \label{eq:14}
  \zeta_X (s,1) = (q-1) q^{-s (1-g)} \zeta_X (s) \; .
\end{equation}

The rest of this section contains observations of a tentative nature which are not necessary for the sequel. It is unknown whether $Z_X (T,u)$ has a natural cohomological interpretation. The properties of $Z_X (T,u)$ are compatible with the following conjectural setup. Let $K$ be a field of characteristic zero containing $\Q (u)$. There may exist a cohomology theory $QH^i$ for varieties over finite fields consisting of finite-dimensional $K$-vector spaces on which the $q$-linear Frobenius $\Fr_q$ induces a $K$-linear map $\Fr^*_q$ such that we have:
\begin{equation}
  \label{eq:15}
  Z_X (T,u) = \prod^2_{i=0} \ddet_K (1 - T \Fr^*_q \tei Q H^i (X))^{(-1)^{i+1}} \; .
\end{equation}
We get the correct denominator in (\ref{eq:4}) if
\[
QH^0 (X) = K \; , \; \Fr^*_q = \id \quad \mbox{and} \quad QH^2 (X) \cong K \; \mbox{with} \; \Fr^*_q = u \cdot \id
\]
and
\[
QH^i (X) = 0 \quad \mbox{for} \; i > 2 \; .
\]
Then $P (T,u)$ would be the characteristic polynomial of $\Fr^*_q$ on $QH^1 (X)$ and therefore we would have
\[
\dim_K QH^1 (X) = 2g \; .
\]
The functional equation (\ref{eq:7}) would be a consequence of Poincar\'e duality -- a perfect $\Fr^*_q$-equivariant pairing of $K$-vector spaces:
\[
QH^i (X) \times QH^{2-i} (X) \longrightarrow QH^2 (X) \cong K \; .
\]
Comparing the logarithmic derivatives of (\ref{eq:3}) and (\ref{eq:15}) at $T = 0$ gives
\begin{equation}
  \label{eq:16}
  \sum^2_{i=0} (-1)^i \Tr (\Fr^*_q \tei QH^i (X)) = \betr{X(\F_q)} \; .
\end{equation}
Moreover Poincar\'e duality would imply
\[
\det (\Fr^*_q \tei QH^1 (X)) = u^g \; .
\]
However, if in (\ref{eq:16}) we replace $\Fr^*_q$ by its power $\Fr^{\nu *}_q$ we do not obtain $\betr{X (\F_{q^{\nu}})}$ for $\nu \ge 2$ if $g \ge 1$. 

\section{Background on two-variable zeta functions of number fields}
\label{sec:3}

We begin by collecting some notions from one-dimensional Arakelov theory following \cite{GS}. \\
For a number field $k / \Q$ let $\eo_k$ be its ring of integers. By $\ep$ we denote the prime ideals in $\eo_k$ and by $v$ the infinite places of $k$. Consider the ``arithmetic curve''
\[
X_k = \spec \eo_k \cup \{ v \tei \infty \} \; .
\]
The elements of the group
\[
Z^1 (X_k) = \bigoplus_{\ep} \Z \cdot \ep \oplus \bigoplus_{v \tei \infty} \R \cdot v
\]
are called Arakelov divisors. Define a map
\[
\div : k^* \longrightarrow Z^1 (X_k)
\]
by the formula 
\[
\div (f) = \sum_{\ep} \ord_{\ep} (f)\ep - \sum_v e_v \log \betr{f}_v  v \; .
\]
Here $\betr{f}_v = \betr{\sigma_v (f)}$ for any embedding $\sigma_v$ in the class $v$ and $e_v = 1$ if $v$ is real and $e_v = 2$ if $v$ is complex. 
The cokernel of $\div$ is called the Arakelov Chow group of $X_k$.

With the evident topologies the groups $k^* , Z^1 (X_k)$ and $CH^1 (X_k)$ become locally compact topological groups. The counting measure on $\bigoplus_{\ep} \Z \cdot \ep$ and the Lebesgue measure on $\bigoplus_{v \tei \infty} \R \cdot v$ induce Haar measures $d D$ on $Z^1 (X_k)$ and $d \Dh$ on $CH^1 (X_k)$.

For an Arakelov divisor
\[
D = \sum_{\ep} \nu_{\ep} \cdot \ep + \sum_v x_v \cdot v \quad \mbox{in} \; Z^1 (X_k)
\]
set
\[
I (D) = \prod_{\ep} \ep^{-\nu_{\ep}} \; .
\]
The infinite components of $D$ determine a norm $\|\;\|_D$ on $k \otimes \R = \bigoplus_v k_v$ by the formula
\[
\| (z_v) \|^2_D = \sum_v \betr{z_v}^2 \Betr{{1}}^2_v \; .
\]
Here $\Betr{{1}}^2_v = e^{-2x_v}$ if $v$ is real and $\Betr{{1}}^2_v = 2e^{-x_v}$ if $v$ is complex.

For $f \in k \hookrightarrow k \otimes \R$ we then have
\begin{equation}
  \label{eq:17}
  \Betr{{f}}^2_D = \sum_{v\,\mathrm{real}} \betr{f}^2_v e^{-2x_v} + 2 \sum_{v\,\mathrm{complex}} \betr{f}^2_v e^{-x_v} \; .
\end{equation}
The embedding $I (D) \hookrightarrow k \otimes \R$ and the norm $\Betr{{\;}}_D$ turn $I (D)$ into a lattice. The lattices $I (D)$ and $I (D')$ are isometric (by an $\eo_k$-linear isometry) if and only if $[D] = [D']$ in $CH^1 (X_k)$. 

Let $\kappa$ be the Arakelov divisor with zeroes at the infinite components and $I (\kappa) = \ed^{-1}$, where $\ed = \ed_{k / \Q}$ is the different of $k / \Q$. 

In the number field case, van der Geer and Schoof replace the order $q^{h^i (\Dh)}$ of $H^i (X , \Oh (D))$ for $X / \F_q$ by the Theta series:
\begin{equation}
  \label{eq:18}
  k^0 (\Dh) = \sum_{f \in I (D)} e^{-\pi \Betr{{f}}^2_D}
\end{equation}
and
\begin{equation}
  \label{eq:19}
  k^1 (\Dh) = k^0 ([\kappa] - \Dh)
\end{equation}
for $\Dh = [D]$ in $CH^1 (X_k)$. For quadratic number fields the behaviour of $k^0 (\Dh)$ is studied in some detail in \cite{F}.

According to \cite{GS} proposition 1, the Poisson summation formula gives the Riemann--Roch type formula
\begin{equation}
  \label{eq:20}
  k^0 (\Dh) k^1 (\Dh)^{-1} = \Nh (\Dh) d^{-1/2}_k \; .
\end{equation}
Here $d_k = \betr{d_{k / \Q}}$ is the absolute value of the discriminant of $k / \Q$ and
\[
\Nh : CH^1 (X_k) \longrightarrow \R^*_+
\]
is the Arakelov norm induced by the map
\[
\Nh : Z^1 (X_k) \longrightarrow \R^*_+ \; , \; \Nh (D) = \prod_{\ep} N\ep^{\nu_{\ep}} \prod_v e^{x_v} \; .
\]
Let $Z^1 (X_k)^0$ be the kernel of this map and set
\[
CH^1 (X_k)^0 = Z^1 (X_k)^0 / \div (k^*) \; .
\]
This is a compact topological group which fits into the exact sequence
\begin{equation}
  \label{eq:21}
  0 \longrightarrow CH^1 (X_k)^0 \longrightarrow CH^1 (X_k) \xrightarrow{\Nh} \R^*_+ \longrightarrow 1 \; .
\end{equation}
Let $d^0 \Dh$ be the Haar measure on $CH^1 (X_k)^0$ with
\begin{equation}
  \label{eq:22}
  \vol (CH^1 (X_k)^0) = hR_k
\end{equation}
where $h = \betr{CH^1 (\spec \eo_k)}$ is the class number of $k$ and $R_k$ is the regulator. Then we have
\begin{equation}
  \label{eq:23}
  d\Dh = d^0 \Dh \frac{dt}{t} \; .
\end{equation}
For $t$ in $\R^*_+$ consider the Arakelov divisor, where $n = (k : \Q)$
\[
D_t = n^{-1} \sum_{v\,\mathrm{real}} \log t \cdot v + n^{-1} \sum_{v\,\mathrm{complex}} 2 \log t \cdot v \; .
\]
Setting $\Dh_t = [D_t]$ we have $\Nh (\Dh_t) = t$, so that the homomorphism $t \mapsto \Dh_t$ provides a splitting of (\ref{eq:21}). 

We need the following estimates:

\begin{prop}
  \label{t31}
For every number field $k$ and every $R \ge 0$ there are positive constants $c_1 , c_2 , \alpha$ such that uniformly in $\Dh \in CH^1 (X_k)^0$ and $\betr{w} \le R$ we have the estimates\\
a) $\betr{k^0 (\Dh + \Dh_t)^w - 1} \le c_1 \betr{w} \exp (- \pi nt^{-2/n})$ \quad for all $0 < t \le \sqrt{d_k}$.\\
b) $\betr{k^0 (\Dh + \Dh_t)^w - t^w d^{-w/2}_k} \le c_2 \betr{w} \exp (- \alpha t^{2/n})$ \quad for all $t \ge \sqrt{d_k}$. 
\end{prop}

\begin{proof}
  According to \cite{GS} corollary 1 there is a constant $\beta > 0$ depending only on the field $k$ such that for all $\Dh$ in $CH^1 (X_k)^0$ and all $0 < t \le d^{1/2}_k$ we have
  \begin{equation}
    \label{eq:24}
    0 < k^0 (\Dh + \Dh_t) - 1 \le \beta \exp (- \pi n t^{-2/n}) \; .
  \end{equation}
We may assume that $R \ge 1$. For every $- \halb \le x \le \halb$ and $\betr{w} \le R$ setting
\begin{equation}
  \label{eq:25}
  (1 + x)^w = 1 + wx + wx^2 \vartheta (x,w)
\end{equation}
we have
\begin{equation}
  \label{eq:26}
  \betr{\vartheta (x,w)} \le e^{2R} \; .
\end{equation}
Namely, writing
\[
(1 + x)^w = e^{w \log (1 + x)} = e^{wx (1 + \eta x)}
\]
we have $\eta = - \halb + \frac{x}{3} - \frac{x^2}{4} + - \ldots$ and hence $\betr{\eta} \le 1$. Expanding $e^{wx (1 + \eta x)}$ as a Taylor series and estimating gives inequality (\ref{eq:26}). For the moment we only need the following consequence of (\ref{eq:26}):
\begin{equation}
  \label{eq:27}
  \betr{(1 + x)^w - 1} \le x \betr{w} \Big( 1 + \halb e^{2R} \Big) \quad \mbox{for} \; 0 \le x \le 1/2 \; \mbox{and} \; \betr{w} \le R \ge 1 \; .
\end{equation}
If $\varepsilon = \varepsilon (k) > 0$ is sufficiently small, (\ref{eq:24}) implies that
\[
x = k^0 (\Dh + \Dh_t) - 1
\]
lies in $( 0, 1/2)$ for all $0 < t \le \varepsilon$ and all $\Dh$. Using (\ref{eq:24}) and (\ref{eq:27}) we therefore find a constant $c'_1$ such that a) holds for all $0 < t \le \varepsilon$. By compactness of
\[
CH^1 (X_k)^0 \times \{ \betr{w} \le R \} \times [\varepsilon , \sqrt{d_k}] 
\]
and continuity of $\frac{1}{w} (k^0 (\Dh + \Dh_t)^w - 1)$ as a function of $\Dh, w$ and $t$ there is a constant $c''_1$ such that a) holds in $\varepsilon \le t \le \sqrt{d_k}$. Thus we get the estimate a) by taking $c_1 = \max (c'_1 , c''_1)$. The estimate b) follows from a) using the Riemann--Roch formula (\ref{eq:20}) and observing that $\Nh ([\kappa]) = d_k$. 
\end{proof}
\bigskip

The two-variable zeta function of van der Geer and Schoof is defined by an integral analogous to the series (\ref{eq:9})
\begin{equation}
  \label{eq:28}
  \zeta^{GS}_{X_k} (s,t) = \int_{CH^1 (X_k)} k^0 (\Dh)^s k^1 (\Dh)^t \, d\Dh \quad \mbox{in} \; \RRe s < 0 , \RRe t < 0 \; .
\end{equation}
According to \cite{GS} proposition 6, this integral defines a holomorphic function in $\RRe s < 0 , \RRe t < 0$. This also follows from the considerations below.

Making the substitution $\Dh \mapsto [\kappa] - \Dh$ in the integral we find the formula
\begin{equation}
  \label{eq:29}
  \zeta^{GS}_{X_k} (s,t) = \int_{CH^1 (X_k)} k^0 (\Dh)^t k^1 (\Dh)^s \, d\Dh \quad \mbox{in} \; \RRe s < 0 , \RRe t < 0 \; .
\end{equation}
We will use the Lagarias--Rains variables $s$ and $w = t+s$ and concentrate on the function
\begin{eqnarray}
  \label{eq:30}
  \zeta_{X_k} (s,w) & = & \zeta^{GS}_{X_k} (s , w-s) = \int_{CH^1 (X_k)} k^0 (\Dh)^{w-s} k^1 (\Dh)^s \, d\Dh \\
& \overset{(\ref{eq:20})}{=} & d^{s/2}_k \int_{CH^1 (X_k)} k^0 (\Dh)^w \Nh (\Dh)^{-s} \, d\Dh \; . \label{eq:31}
\end{eqnarray}
It is holomorphic in the region $\RRe w < \RRe s < 0$.

Most of the following proposition is stated in \cite{GS} and proved in \cite{LR} Appendix using references to Ch. XIII of Serge Lang's book on algebraic number theory. Below we will write down the direct proof which is implicit in \cite{GS}.

\begin{prop}
  \label{t32}
  The function $\zeta_{X_k} (s,w)$ has a meromorphic continuation to $\C^2$ and it satisfies the functional equation
\[
\zeta_{X_k} (s,w) = \zeta_{X_k} (w-s,w) \; .
\]
Moreover the function
\[
w^{-1} s (w-s) \zeta_{X_k} (s,w)
\]
is holomorphic in $\C^2$. More precisely, the integral
\[
J (s,w) = \int^{\sqrt{d_k}}_0 \int_{CH^1 (X_k)^0} w^{-1} (k^0 (\Dh + \Dh_t)^w - 1) d^0 \Dh \; t^{-s} \frac{dt}{t}
\]
defines an entire function in $\C^2$ and we have the formula
\[
\zeta_{X_k} (s,w) = w \left( d^{s/2}_k J (s,w) + d^{(w-s)/2}_k J (w-s,w) \right) - \left( \frac{1}{s} + \frac{1}{w-s} \right) hR_k \; .
\]
Recall that $\vol CH^1 (X_k)^0 = hR_k$. Finally, for $w = 1$ one has
\begin{equation}
  \label{eq:32}
  \zeta_{X_k} (s,1) = \betr{\mu (k)} d^{s/2}_k 2^{-r_1/2} \hat{\zeta}_k (s) \; .
\end{equation}
\end{prop}

Here $\hat{\zeta}_k (s)$ is the completed Dedekind zeta function of $k$
\[
\hat{\zeta}_k (s) = \zeta_k (s) \Gamma_{\R} (s)^{r_1} \Gamma_{\C} (s)^{r_2}
\]
where we have set
\[
\Gamma_{\R} (s) = 2^{-1/2} \pi^{-s/2} \Gamma (s/2) \quad \mbox{and} \quad \Gamma_{\C} (s) = (2 \pi)^{-s} \Gamma (s) \; .
\]
Thus $\Gamma_{\R} (s) \Gamma_{\R} (s+1) = \Gamma_{\C} (s)$. Here $r_1$ and $r_2$ are the numbers of real resp. complex places of $k$. Moreover $\mu (k)$ is the group of roots of unity in $k$.

\begin{rems}
  {\bf 1} Formula (\ref{eq:32}) coincides with the corresponding formula in \cite{GS} proposition 6 after correcting two small misprints in that paper: We have $\sqrt{\betr{\Delta}}^s$ instead of $\sqrt{\betr{\Delta}}^{s/2}$ in \cite{GS} proposition 6 and $2^{-1} \pi^{-s/2} \ldots$ instead of $2 \pi^{-s/2} \ldots$ in the third equality on p. 388 above of \cite{GS}.\\
{\bf 2} The reason for our normalization of $\Gamma_{\R} (s)$ comes from the theory of zeta-regularization, c.f. section 5.
\end{rems}

\begin{proof}
  We write the integral representation (\ref{eq:31}) for $\zeta_{X_k} (s,w)$ as a sum of two contributions:
  \begin{equation}
    \label{eq:33}
    \zeta_{X_k} (s,w) = I (s,w) + II (s,w)
  \end{equation}
where 
\[
I (s,w) = d^{s/2}_k \int^{\sqrt{d_k}}_0 \int_{CH^1 (X_k)^0} k^0 (\Dh + \Dh_t)^w d^0 \Dh \; t^{-s} \frac{dt}{t}
\]
and
\[
II (s,w) = d^{s/2}_k \int^{\infty}_{\sqrt{d_k}} \int_{CH^1 (X_k)^0} k^0 (\Dh + \Dh_t)^w d^0 \Dh \; t^{-s} \frac{dt}{t} \; .
\]
The estimate in proposition \ref{t31} a) shows that the first integral defines a holomorphic function in $\RRe s< 0 , w \in \C$. Here and in the following we use the following well known fact. Consider a function $f (s,x)$ holomorphic in $s$ and $\mu$-integrable in $x$ which locally in $s$ is bounded by integrable functions of $x$. Then the integral $\int f (s,x) d \mu (x)$ is holomorphic in $s$. Writing $I (s,w)$ in the form 
\begin{equation}
  \label{eq:34}
\small  I (s,w) = d^{s/2}_k \int^{\sqrt{d_k}}_0 \int_{CH^1 (X_k)^0} (k^0 (\Dh + \Dh_t)^w - 1) d^0 \Dh \; t^{-s} \frac{dt}{t}  - \frac{hR_k}{s}
\end{equation}
the same estimate gives its meromorphic continuation to $\C^2$. Note that, even divided by $w$ the first term is holomorphic in $\C^2$. 

Using Riemann--Roch (\ref{eq:20}) a short calculation shows that for $\RRe s > \RRe w$ we have
\begin{equation}
  \label{eq:35}
  II (s,w) = I (w-s,w) \; .
\end{equation}
In particular the integral (\ref{eq:31}) defines a holomorphic function in $\RRe w < \RRe s < 0$ as asserted earlier. Using (\ref{eq:34}) we find the formula:
\begin{equation}
  \label{eq:36}
\small  II (s,w) = d^{(w-s)/ 2}_k \int^{\sqrt{d_k}}_0 \int_{CH^1 (X_k)^0} (k^0 (\Dh + \Dh_t)^w - 1) d^0 \Dh \; t^{- (w-s)} \frac{dt}{t} - \frac{hR_k}{w-s}
\end{equation}
which gives the meromorphic continuation of $II (s,w)$ to $\C^2$: Again, even after division by $w$ the first term is holomorphic in $\C^2$. This implies the assertions of the proposition except for formula (\ref{eq:32}) which requires a lemma that will be useful in the next section as well:
\end{proof}

\begin{lemma}
  \label{t33}
In the region $\RRe s > \RRe w , \RRe s > 0$ the following integral representation holds, the integral defining a holomorphic function even after division by $w$:
\begin{equation}
  \label{eq:37}
  \zeta_{X_k} (s,w) = d^{s/2}_k \int_{CH^1 (X_k)} (k^0 (\Dh)^w - 1) \Nh \Dh^{-s} d \Dh \; .
\end{equation}
\end{lemma}

\begin{proofof}
  {\bf formula (\ref{eq:32})} Using (\ref{eq:37}) we find for $w = 1 < \RRe s$ that
\[
\betr{\mu (k)}^{-1} d^{-s/2}_k \zeta_{X_k} (s,1) = \betr{\mu (k)}^{-1} \int_{CH^1 (X_k)} (k^0 (\Dh) - 1) \Nh \Dh^{-s} d \Dh \; .
\]
Now on p. 388 of \cite{GS} this integral is shown to equal
\[
(2^{-1} \pi^{-s/2} \Gamma (s/2))^{r_1} ((2 \pi)^{-s} \Gamma (s))^{r_2} \zeta_k (s)
\]
c.f. remark 1 above.
\end{proofof}

\begin{proofof}
  {\bf the lemma} The estimate in proposition \ref{t31}, b) shows that the following formula is valid in the region $\RRe s > \RRe w , \RRe s > 0$:
  \begin{equation}
    \label{eq:38}
\small    II (s,w) = d^{s/2}_k \int^{\infty}_{\sqrt{d_k}} \int_{CH^1 (X_k)^0} (k^0 (\Dh + \Dh_t)^w-1) d^0 \Dh \, t^{-s} \frac{dt}{t} + \frac{hR_k}{s} \; .
  \end{equation}
The integral defines a holomorphic function in this region even after division by $w$. As the integral in formula (\ref{eq:34}) for $w^{-1} I (s,w)$ gives a holomorphic function in $\C^2$ the assertion follows by adding equations (\ref{eq:34}) and (\ref{eq:38}).
\end{proofof}

\begin{rem}
  For $k = \Q$ a more elaborate version of the lemma is given in \cite{LR} Theorem 2.2.
\end{rem}

Proposition \ref{t32} and formula (\ref{eq:32}) in particular suggest that a better definition of a two variable zeta function might be the following
\[
  \zeta (X_k , s,w) = w^{-1} \frac{2^{r_1/2}}{\betr{\mu (k)}} d^{-s/2}_k \zeta_{X_k} (s,w) \; .
\]
This is a meromorphic function on $\C^2$ which satisfies the equations
\begin{equation}
\label{eq:39}
\zeta (X_k , w-s , w) = d^{s-w/2}_k \zeta (X_k , s,w) \quad \mbox{and} \quad
\zeta (X_k ,s,1) = \hat{\zeta}_k (s) \; .
\end{equation}
In section \ref{sec:5} we will see that $\zeta (X_k ,s,w)$ is the ``$\frac{1}{2\pi}$-zeta regularized version'' of $\zeta_{X_k} (s,w)$. We also consider an entire version of this function which in the one variable case and in \cite{LR} is called the $\xi$-function. Because of our different normalization we give it another name which is suggested by the cohomological arguments in section \ref{sec:6}. 

\begin{defn}
  \label{t34}
The two-variable $L$-function of $X_k$ is defined by the formula
\begin{eqnarray*}
  L (H^1 (X_k) , s,w) & = & \frac{s}{2 \pi} \frac{s-w}{2\pi} \zeta (X_k ,s,w)\\
& = & \frac{1}{4\pi^2} \frac{s (s-w)}{w} \frac{2^{r_1/2}}{\betr{\mu (k)}} d^{-s/2}_k \zeta_{X_k} (s,w) \; .
\end{eqnarray*}
According to proposition \ref{t32} it is holomorphic in $\C^2$ and satisfies the functional equation 
\[
L (H^1 (X_k) , w-s , w) = d^{s-w/2}_k L (H^1 (X_k) , s,w) \; .
\]
\end{defn}

\begin{prop}
  \label{t35}
For any $k / \Q$ and every fixed $w$ the entire function \\
$L (H^1 (X_k) , s,w)$ of $s$ has order at most one.
\end{prop}

\begin{proof}
  Proposition \ref{t32} implies the formula
\[
\small L (H^1 (X_k) , s,w) = s (s-w) (T (s,w) + d^{\frac{w}{2}-s}_k T (w-s ,w )) + \frac{d^{-s/2}_k}{4\pi^2} \frac{2^{r_1/2}}{\betr{\mu (k)}} hR_k
\]
where $T (s,w)$ is the entire function in $\C^2$ defined by the integral
\[
T (s,w) = \frac{1}{4\pi^2} \frac{2^{r_1/2}}{\betr{\mu (k)}} \int^{\sqrt{d_k}}_0 \int_{CH^1 (X_k)^0} w^{-1} (k^0 (\Dh + \Dh_t)^w-1) d^0 \Dh \, t^{-s} \frac{dt}{t} \; .
\]
Using the estimate in proposition \ref{t31}, a) we find for some $c (w)> 0$:
\begin{eqnarray*}
  \betr{T (s,w)} & \le & c (w) \int^{\sqrt{d}_k}_0 \exp (-\pi n t^{-2/n}) t^{-\RRe s} \frac{dt}{t} \\
& = & c (w) d^{-\RRe s /2}_k \int^{\infty}_1 \exp ( -\pi n d^{-1/n}_k t^{2/n}) t^{\RRe s} \frac{dt}{t} \; .
\end{eqnarray*}
For $\RRe s \le 1$ the latter integral is bounded. For $\RRe s > 1$ we have
\begin{eqnarray*}
  \betr{T (s,w)} & \le & c (w) d^{-\RRe s/2}_k \int^{\infty}_0 \exp (-\pi n d^{-1/n}_k t^{2/n}) t^{\RRe s} \frac{dt}{t} \\
& = & \frac{nc (w)}{2} (\pi n)^{- \frac{n\RRe s}{2}} \Gamma \left( \frac{n \RRe s}{2} \right) \\
& = & O \left( \exp \left( \frac{n}{2} \RRe s\right) \log (\RRe s) \right)
\end{eqnarray*}
where the $O$-constant depends on $w$. Hence for all $s \in \C$ we have
\[
\betr{T (s,w)} = O \left( \exp \left( \frac{n}{2} \betr{s} \log \betr{s} \right) \right) \; .
\]
Thus for every $\varepsilon > 0$ the required estimate holds:
\[
\betr{L (H^1 (X_k) , s,w)} = O (\exp (\betr{s}^{1 + \varepsilon})) \quad \mbox{for} \; s \in \C \; .
\]
\end{proof}

\begin{rem}
  For $k = \Q$ Lagarias and Rains prove that $L (H^1 (X_{\Q}) ,s,w)$ is entire of order at most one as a function of two variables, \cite{LR} Theorem 4.1. They also mention that this assertion holds for general $k$ as well.
\end{rem}

\section{An oscillatory integral in the geometry of numbers}
\label{sec:4}

Recall that an Arakelov divisor $D$ in $Z^1 (X_k)$ may be viewed as the lattice $(I (D) , \Betr{{\;}}_D)$. Two divisors define the same class $\Dh$ in $CH^1 (X_k)$ if and only if the corresponding lattices are isometric. In particular the following numbers are well defined for $\Dh = [D]$:
\begin{eqnarray*}
  a (\Dh) & = & \min \{ \Betr{{f}}^2_D \tei 0 \neq f \in I (D) \} \\
b (\Dh) & = & \min \{ \Betr{{f}}^2_D \tei f \in I (D) \; \mbox{such that} \; \Betr{{f}}^2_D > a (\Dh) \} \\
\nu (\Dh) & = & \betr{\{ f \in I (D) \tei \Betr{{f}}^2_D = a_{\Dh}\} } \; .
\end{eqnarray*}
By definition $b (\Dh) > a (\Dh) > 0$ are positive real numbers and $\nu (\Dh)$ is a positive integer -- the so called kissing number of the lattice class.

These numbers arise naturally in the study of theta functions: Ordering terms, we may write
\begin{eqnarray*}
  k^0 (\Dh + \Dh_t) & = & \sum_{f \in I (D)} \exp (-\pi t^{-2/n} \Betr{{f}}^2_D) \\
& = & 1 + \nu (\Dh) e^{-\pi t^{-2/n}} a (\Dh) + \ldots
\end{eqnarray*}
Here the next term is $e^{-\pi t^{-2/n} b (\Dh)}$ with its multiplicity.

\begin{prop}
  \label{t41}
On $CH^1 (X_k)$ the function $a$ is continuous whereas $b$ and $\nu$ are only upper semicontinuous. In particular $a,b$ and $\nu$ are measurable. We have $b (\Dh) \le 4a (\Dh)$ for all $\Dh$, and $\nu$ is locally bounded. On $CH^1 (X_k)^0$ the functions $a,b, \nu$ are bounded.
\end{prop}

Points of discontinuity for $b$ and $\nu$ arise as follows. Already for $k = \Q (\sqrt{2})$ there exist convergent sequences $\Dh_n \to \Dh$ even in $CH^1 (X_k)^0$ such that $b (\Dh_n) \to a (\Dh)$. Thus at the point $\Dh$ we have $\lim_{n\to\infty} b (\Dh_n) < b (\Dh)$ and also the multiplicity $\nu$ jumps up. 

\begin{proof}
  Fix an element $f \in I (D)$ with $\Betr{{f}}^2_D = a_{\Dh}$. Then $\Betr{{2f}}^2_D = 4 a_{\Dh}$. Thus $b_{\Dh} \le 4a_{\Dh}$. The continuity properties may be checked locally. So let us fix a class $\Dh^0 = [D^0]$ in $CH^1 (X_k)$ and write:
\[
D^0 = \sum_{\ep} \nu^0_{\ep} \cdot \ep + \sum_v x^0_v \cdot v \quad \mbox{in} \; Z^1 (X_k) \; .
\]
Let $V$ be an open neighborhood of $x^0 = (x^0_v)_{v \tei \infty}$ in $\bigoplus_{v \tei \infty} \R$ and consider the continuous map:
\[
V \longrightarrow CH^1 (X_k) \; , \; x \longmapsto \Dh_x = [D_x] \quad \mbox{where} \; D_x = \sum_{\ep} \nu^0_{\ep} \cdot \ep + \sum_{v \tei \infty} x_v \cdot v \; .
\]
For $V$ small enough this map is a homeomorphism of $V$ onto an open neighborhood $U$ of $\Dh^0$ in $CH^1 (X_k)$. Fix some $R > 0$ such that for all $x$ in $V$ we have
\[
R^{-1} \le e^{x_v} \le R \; \mbox{if $v$ is real \quad and} \quad R^{-2} \le e^{x_v} \le R^2 \; \mbox{if $v$ is complex} \; .
\]
It follows that for $x \in V$ and all $f \in I (D_x) = I (D^0)$ we have the estimate
\begin{equation}
  \label{eq:40}
\small  R^{-2} \Betr{f}^2 \le \Betr{{f}}^2_{D_x} = \sum_{v \, \mathrm{real}} \betr{f}^2_v e^{-2x_v} + 2 \sum_{v\, \mathrm{complex}} \betr{f}^2_v e^{-x_v} \le 2R^2 \Betr{f}^2 \; .
\end{equation}
Here
\[
\Betr{{f}} = \Big( \sum_{v \tei \infty} \betr{f}^2_v \Big)^{1/2}
\]
is the Euklidean norm in $k \otimes \R$ applied to the element $f \in k \subset k \otimes \R$.

Since $I (D^0)$ is discrete in $k \otimes \R$ it follows that for any $C > 0$ the set
\[
\Fh_C = \{ f \in I (D^0) \tei 0 < \Betr{{f}}^2_{D_x} \le C \quad \mbox{for some} \; x \in V \}
\]
is finite. If $V$ is bounded it also follows that the map $\Dh \mapsto a (\Dh)$ is bounded on $U$ and so is $b$ since $b (\Dh) \le 4a (\Dh)$. Thus for large enough $C > 0$ the finite subset $\Fh = \Fh_C \subset I (D)$ has the following properties: For all $x \in V$ we have:
\begin{eqnarray}
a (\Dh_x) & = & \min \{ \Betr{{f}}^2_{D_x} \tei f \in \Fh \} \nonumber \\
b (\Dh_x) & = & \min \{ \Betr{{f}}^2_{D_x} \tei f \in \Fh \; \mbox{such that} \; \Betr{{f}}^2_{D_x} > a (\Dh_x) \} \nonumber \\
\nu (\Dh_x) & = & \betr{ \{ f \in \Fh \tei \Betr{{f}}^2_{D_x} = a (\Dh_x) \} } \; . \label{eq:41}
\end{eqnarray}
The functions $x \mapsto \Betr{{f}}^2_{D_x}$ for $f \in \Fh$ being continuous it is now clear that $a (\Dh)$ is continuous near $\Dh^0$, hence everywhere since $\Dh^0$ was arbitrary. (This fact is already mentioned in \cite{GS}.)

To check upper semicontinuity of $b$ and $\nu$ at $\Dh^0$ let $\Fh'$ be the subset of $\Fh$ consisting of all $f$ with $\Betr{{f}}^2_{D^0} > a (\Dh_0)$. For small enough $V$ we then have
\begin{equation}
  \label{eq:42}
  \Betr{{f}}^2_{D_x} > a (\Dh_x) \quad \mbox{for all} \; f \in \Fh' \; \mbox{and} \; x \in V
\end{equation}
since both sides are continuous in $x$. It follows that
\[
b (\Dh_x) \le \min \{ \Betr{{f}}^2_{D_x} \tei f \in \Fh' \} =: \mu (x) \; .
\]
Then $\mu$ is continuous and $\mu (x^0) = b (\Dh^0)$. Hence, for every $\varepsilon > 0$ there exists an open neighborhood $V'$ of $x^0$ in $V$ such that
\[
\mu (V') \subset (\mu (x^0) - \varepsilon , \mu (x^0) + \varepsilon) = (b (\Dh^0) - \varepsilon , b (\Dh^0) + \varepsilon) \; .
\]
Thus $b (\Dh_x) \le b (\Dh^0) + \varepsilon$ for all $x \in V'$ and hence $b (\Dh) \le b (\Dh^0) + \varepsilon$ for all $\Dh$ in a neighborhood (the image of $V'$) of $\Dh^0$ in $CH^1 (X_k)$. Hence $b$ is upper semicontinuous at $\Dh^0$.

As for $\nu$, the representation (\ref{eq:41}) shows that $\nu (\Dh_x) \le \betr{\Fh}$ for all $x \in V$. Hence $\nu$ is a locally bounded function on $CH^1 (X_k)$.

With notations as above we have by (\ref{eq:41}) that
\[
\nu (\Dh_x) \le \betr{\Fh \ohne \Fh'} = \nu (\Dh^0) \quad \mbox{for all} \; x \in V \; .
\]
This implies that $\nu$ is upper semicontinuous at $\Dh^0$.
\end{proof}

The following theorem shows that on $CH^1 (X_k)^0$ the function $a = a (\Dh)$ acquires a unique global minimum at $\Dh = 0$. We also describe $a (\Dh)$ explicitly in a neighborhood of $\Dh = 0$. 

Set
\[
a_{\min} = \min \{ a (\Dh) \tei \Dh \in CH^1 (X_k)^0 \} > 0 
\]
and
\[
b_{\inf} = \inf \{ b (\Dh) \tei \Dh \in CH^1 (X_k)^0 \} \; .
\]

\begin{theorem}
  \label{t42}
Set $n = (k : \Q)$ and let the notations be as above. \\
{\bf 1} $a_{\min} = n$.\\
{\bf 2} For $\Dh \in CH^1 (X_k)^0$ we have $a (\Dh) = a_{\min}$ if and only if $\Dh = 0$.\\
{\bf 3} For the representative $D = 0$ of $\Dh = 0$ and $f \in \eo_k = I (0)$ we have $\Betr{{f}}^2_0 = a (0) = a_{\min}$ if and only if $f \in \mu (k)$.\\
{\bf 4} For $\Dh \in CH^1 (X_k)^0$ with $I (D)$ non-principal there is the estimate
\[
a (\Dh) \ge \sqrt[n]{4} \; a_{\min} = n \sqrt[n]{4} \; .
\]
{\bf 5} For every open neighborhood $U$ of $\Dh = 0$ in $CH^1 (X_k)^0$ there is a positive $\varepsilon$ such that $a (\Dh) < a_{\min} + \varepsilon$ for some $\Dh \in CH^1 (X_k)^0$ implies $\Dh \in U$.\\
{\bf 6} There is a neighborhood $U$ of $\Dh = 0$ in $CH^1 (X_k)^0$ with the following properties: Every $\Dh \in U$ has the form $\Dh = [D]$ with $D = \sum_{v \tei \infty} x_v \cdot v$. For $f \in I (D) = \eo_k$ we have:
\[
\Betr{{f}}^2_D = a (\Dh) \quad \mbox{if and only if} \; f \in \mu (k) \; .
\]
Moreover:
\[
a (\Dh) = \sum_{v\,\real} e^{-2x_v} + 2 \sum_{v\,\complex} e^{-x_v}
\]
and $\nu (\Dh) = \betr{\mu (k)}$.\\
{\bf 7} We have $b_{\inf} > a_{\min}$. 
\end{theorem}

\begin{proof}
  The main tool is the inequality between the arithmetic and the geometric mean. This inequality was already used in \cite{GS}. Let $\Betr{{\;}}_v = \betr{\;}^{e_v}_v$ be the normalized absolute value at the infinite place $v$.\\
{\bf 1} For $\Dh = [D]$ in $CH^1 (X_k)^0$ and $f \in I (D)$ we have
\begin{eqnarray*}
  \Betr{{f}}^2_D & = & \sum_{v \, \real} (\Betr{{f}}_v e^{-x_v})^2 + \sum_{v\,\complex} \Betr{{f}}_v e^{-x_v} + \sum_{v\,\complex} \Betr{{f}}_v e^{-x_v} \\
& \overset{\mathrm (a)}{\ge} & n \left( \prod_v \Betr{{f}}_v \right)^{2/n} \left( \prod_v e^{-x_v} \right)^{2/n} = n \betr{N (f)}^{2/n} \left( \prod_v e^{x_v} \right)^{-2/n} \\
& = & n (\betr{N (f)} / N (I (D)))^{2 /n } \; .
\end{eqnarray*}
Here (a) is the arithmetic-geometric mean inequality and we have used that
\[
1 = \Nh (\Dh) = \prod_{\ep} N\ep^{\nu_{\ep}} \prod_v e^{x_v} = N (I (D))^{-1} \prod_v e^{x_v} \; .
\]
Now $I (D)$ divides $(f)$ and for $f \neq 0$ we therefore have
\[
\betr{N (f)} / N (I (D)) \ge 1 \; .
\]
It follows that $\Betr{{f}}^2_D \ge n$, so that $a (\Dh) \ge n$ and therefore $a_{\min} \ge n$. On the other hand for $D = 0$ and $f \in \mu (k)$ we have $\Betr{{f}}^2_0 = r_1 + 2r_2 = n$. Therefore $a (0) = n$ and hence $a_{\min} = n$.\medskip

\noindent {\bf 2} We have seen that $a (0) = a_{\min}$. Now assume that $a (\Dh) = a_{\min}$. Then there is some $f \in I (D)$ with $\Betr{{f}}^2_D = n$. It follows that $\betr{N (f)} = N (I (D))$ hence that $I (D) = (f)$ is principal and that we have equality in (a) above. Now in the arithmetic-geometric mean inequality, equality is achived precisely if all terms are equal. Thus there is a positive real $\xi$ such that
\[
\xi = (\Betr{{f}}_v e^{-x_v})^2 \quad \mbox{for real $v$ and} \; \xi = \Betr{{f}}_v e^{-x_v} \; \mbox{for complex} \; v \; .
\]
Hence
\[
\xi^n = \xi^{r_1} \xi^{2r_2} = \left( \prod_v \Betr{{f}}_v e^{-x_v} \right)^2 = (\betr{N (f)} N (I (D))^{-1})^2 = 1
\]
since $\Nh (\Dh) = 1$ and $\betr{N (f)} = N (I (D))$ as observed above. Thus $\xi = 1$ and therefore
\begin{equation}
  \label{eq:43}
  \Betr{{f}}_v = e^{x_v} \quad \mbox{for all} \; v \tei \infty \; .
\end{equation}
It follows that
\begin{eqnarray*}
  \div f^{-1} & = & \sum_{\ep} \ord_{\ep} f^{-1} \cdot \ep - \sum_v \log \Betr{{f^{-1}}}_v \cdot v \\
& = & \sum_{\ep} \ord_{\ep} I (D)^{-1} \cdot \ep + \sum_v \log \Betr{{f}}_v \cdot v \\
& = & \sum_{\ep} \nu_{\ep} \cdot \ep + \sum_v x_v \cdot v = D \; .
\end{eqnarray*}
Hence $\Dh = [D] = 0$ in $CH^1 (X_k)^0$ and {\bf 2} is proved.
\medskip

\noindent {\bf 3} For $D = 0$ we have $I (D) = \eo_k$. For $f \in I (D) = \eo_k$ the equation $\Betr{{f}}^2_0 = a (0) = n$ implies $\Betr{{f}}_v = 1$ for all $v \tei \infty$ by (\ref{eq:43}). Since $\Betr{{f}}_{\ep} \le 1$ for all finite primes $\ep$ it follows by a theorem of Kronecker that $f$ is a root of unity. \medskip

\noindent {\bf 4} If $I (D)$ is non-principal and $0 \neq f \in I (D)$, then we have $(f) = I (D) \cdot \ea$ for some integral ideal $\ea \neq \eo_k$. Hence $\betr{N (f)} \ge 2 N (I (D))$ and {\bf 4} follows from the above estimate for $\Betr{{f}}^2_D$.\medskip

\noindent {\bf 5} Let $a_U$ be the minimum of the continuous function $a = a (\Dh)$ on the compact set $CH^1 (X_k)^0 \ohne U$. For $\Dh \neq 0$ we have $a (\Dh) > a_{\min}$ by {\bf 2}. Hence $\varepsilon := a_U - a_{\min} > 0$. It is clear that $a (\Dh) < a_{\min} + \varepsilon$ implies that $\Dh \in U$. \medskip

\noindent {\bf 6} As in the proof of proposition \ref{t41} there exists an open neighborhood $V'$ of $x^0 = 0$ in $\{ x \in \bigoplus_{v \tei \infty} \R \tei \sum x_v = 0 \}$ such that firstly the map
\[
V' \longrightarrow CH^1 (X_k)^0 \; , \; x \longmapsto \Dh_x = [D_x] \; \mbox{where} \; D_x = \sum_{v \tei \infty} x_v \cdot v
\]
is a homeomorphism onto an open neighborhood $U'$ of $\Dh = 0$ in $CH^1 (X_k)^0$. In particular $I (D_x) = \eo_k$ for all $x \in V'$. Secondly there is a finite subset $\Fh \supset \mu (k)$ of $\eo_k$ such that for all $x \in V'$ we have:
\begin{eqnarray}
  a (\Dh_x) & = & \min \{ \Betr{{f}}^2_{D_x} \tei f \in \Fh \} \nonumber \\
b (\Dh_x) & = & \min \{ \Betr{{f}}^2_{D_x} \tei f \in \Fh \; \mbox{such that} \; \Betr{{f}}^2_{D_x} > a (\Dh_x) \} \label{eq:44}
\end{eqnarray}
and
\[
\nu (\Dh_x) = \betr{ \{ f \in \Fh \tei \Betr{{f}}^2_{D_x} = a (\Dh_x) \} } \; .
\]
Now, according to {\bf 3} we have
\[
\Betr{{f}}^2_{D_0} = a (\Dh_0) \quad \mbox{for} \; f \in \mu (k)
\]
and
\[
\Betr{{f}}^2_{D_0} > a (\Dh_0) \quad \mbox{for} \; f \in \Fh \ohne \mu (k) \; .
\]
Choose some $\varepsilon > 0$, such that $\Betr{{f}}^2_{D_0} - a (\Dh_0) \ge 2 \varepsilon$ for all $f \in \Fh \ohne \mu (k)$. By a continuity argument we may find an open neighborhood $0 \in V \subset V'$ such that for all $x \in V$ we have
\[
\Betr{{f}}^2_{D_x} - a (\Dh_x) < \varepsilon \quad \mbox{if} \; f \in \mu (k)
\]
and
\[
\Betr{{f}}^2_{D_x} - a (\Dh_x) \ge \varepsilon \quad \mbox{if} \; f \in \Fh \ohne \mu (k) \; .
\]
As $\Betr{{f}}^2_{D_x} = \Betr{{1}}^2_{D_x}$ for all $f \in \mu (k)$ it follows that for $x \in V$ we have
\[
\Betr{{f}}^2_{D_x} = a (\Dh_x) \quad \mbox{if and only if} \; f \in \mu (k) \; .
\]
Moreover $\nu (\Dh_x) = \betr{\mu (k)}$ and
\[
a (\Dh_x) = \Betr{{1}}^2_{D_x} = \sum_{v \, \real} e^{-2x_v} + 2 \sum_{v \, \complex} e^{-x_v} \; .
\]
Therefore, in {\bf 6} we may take $U$ to be the image of $V$ in $CH^1 (X_k)^0$. \medskip

\noindent {\bf 7} Assume that $b_{\inf} = a_{\min}$ and let $(\Dh_n)$ be a sequence of $\Dh_n \in CH^1 (X_k)^0$ with $b (\Dh_n) \to a_{\min}$. Since $CH^1 (X_k)^0$ is compact we may assume that $(\Dh_n)$ is convergent, $\Dh_n \to \Dh_0$. Because of $a_{\min} \le a (\Dh_n) \le b (\Dh_n)$ it follows that $a (\Dh_n) \to a_{\min}$. On the other hand since $a$ is continuous we have $a (\Dh_n) \to a (\Dh_0)$. Hence $a (\Dh_0) = a_{\min}$ and by {\bf 2} this implies that $\Dh_0 = 0$. Thus we have $\Dh_n \to 0$ and $b (\Dh_n) \to a_{\min}$. 

Let $V,U$ and $\Fh$ be as in the proof of {\bf 6}. Then for $f \in \Fh$ and $x \in V$ we have
\[
\Betr{{f}}^2_{D_x} > a (\Dh_x) \quad \mbox{if and only if} \; f \notin \mu (k) \; .
\]
By (\ref{eq:44}) this gives
\[
b (\Dh_x) = \min \{ \Betr{{f}}^2_{D_x} \tei f \in \Fh \ohne \mu (k) \} \quad \mbox{for all} \; x \in V \; .
\]
In particular $b (\Dh_x)$ is a continuous function of $x \in V$ and therefore $b \, |_U$ is continuous. Let $\tU \subset U$ be a compact neighborhood of $\Dh = 0$ in $CH^1 (X_k)^0$. Then there is some $\tilde{\Dh} \in \tU$ with $b (\Dh) \ge b (\tD) > a (\tD) \ge a_{\min}$ for all $\Dh$ in $\tU$. On the other hand, for $n$ large enough we have $\Dh_n \in \tU$ and hence $b (\Dh_n) \ge b (\tD) > a_{\min}$. Hence $b (\Dh_n)$ cannot converge to $a_{\min}$, Contradiction.
\end{proof}

\begin{rem}
Since $\mu (k)$ acts isometrically on $(I (D), \Betr{{\;}}_D)$ and since $\nu (0) = \betr{\mu (k)}$ the minimal value of the function $\nu = \nu (\Dh)$ is $\betr{\mu (k)}$. As $\nu$ is upper semicontinuous it follows that the set of $\Dh$ in $CH^1 (X_k)$ resp. $CH^1 (X_k)^0$ with $\nu (\Dh) = \betr{\mu (k)}$ is open. It should be possible to show that the complements have measure zero.  
\end{rem}

In the following we will deal with the asymptotic behaviour of certain functions defined at least in $\RRe s > 0$ as $\RRe s$ tends to infinity. For such functions $f$ and $g$ we will write
\[
f \sim g \quad \mbox{to signify that} \; \lim_{\RRe s \to \infty} f (s) / g (s) = 1 \; .
\]
The following theorem is the main result of the present section:

\begin{theorem}
  \label{t43}
For a number field $k / \Q$ let $r = r_1 + r_2 - 1$ be the unit rank. Then the entire function
\[
C (s) = \int_{CH^1 (X_k)^0} \nu (\Dh) a (\Dh)^{-s} d^0 \Dh
\]
has the following asymptotic behaviour as $\RRe s \to \infty$
\[
C (s) \sim \betr{\mu (k)} \alpha_k s^{-r/2} n^{-s}  \; .
\]
Here we have set:
\[
\alpha_k = (\pi n)^{r/2} 2^{-r_1 / 2} \sqrt{2/n} \; .
\]
\end{theorem}

\begin{proof}
  If $r = 0$ then $\alpha_k = 1$ and $CH^1 (X_k)^0 = CH^1 (\spec \eo_k)$ is the class group of $k$. Hence $C (s)$ is a finite Dirichlet series. For $k = \Q$ we have $C (s) = \nu (0) a (0)^{-s} = \betr{\mu (\Q)} = 2$. For $k$ imaginary quadratic the main contribution as $\RRe s \to \infty$ comes from the term corresponding to $\Dh = 0$ which is $\nu (0) a (0)^{-s} = \betr{\mu (k)} 2^{-s}$. These assertions follow from theorem \ref{t42} parts {\bf 1} and {\bf 2} (or {\bf 4}) and {\bf 3}.

Now assume that $r \ge 1$. The function $\nu = \nu (\Dh)$ is measurable and bounded on $CH^1 (X_k)^0$ by proposition \ref{t41}. The function $a = a (\Dh)$ is continuous and $CH^1 (X_k)^0$ is compact. Hence $C (s)$ is an entire function of $s$. We will compare $C (s)$ with certain integrals over unbounded domains which can be evaluated explicitely in terms of $\Gamma$-functions. It is not obvious that these integrals converge. For this we require the following lemma where for $x \in \R^N$ we set $\Betr{{x}}_{\infty} = \max \betr{x_i}$. 
\medskip

\noindent After a series of auxiliary results the proof of theorem \ref{t43} is concluded after the proof of corollary \ref{t434} below.

\begin{ulemma}
  \label{t431}
Assume $N \ge 2$ and consider the hyperplane \\
$H_N = \{ x \tei \sum x_i = 0 \}$ in $\R^N$. For every $x$ in $H_N$ we have
\[
\max x_i \ge (N-1)^{-1} \Betr{{x}}_{\infty} \quad \mbox{and} \quad \min x_i \le - (N-1)^{-1} \Betr{{x}}_{\infty} \; .
\]
\end{ulemma}

\begin{proof}
  We may assume that $x_1 \le \ldots \le x_N$, so that $x_1 = \min x_i$ and $x_N = \max x_i$. As $x \in H_N$ we have $x_1 \le 0 \le x_N$. It is clear that $\Betr{{x}}_{\infty} = \max (-x_1 , x_N)$. 

If $\Betr{{x}}_{\infty} = x_N$ the first estimate is clear. If $\Betr{{x}}_{\infty} = -x_1$ then
\[
(N-1) \max x_i = (N-1) x_N \ge x_N + x_{N-1} + \ldots + x_2 = - x_1 = \Betr{{x}}_{\infty} \; .
\]
Hence the first estimate holds in this case as well. The second estimate follows by replacing $x$ with $-x$.
\end{proof}

We can now evaluate a certain class of integrals which are useful for our purposes.

\begin{uprop}
  \label{t432}
For $N \ge 2$ let $d\lambda$ be the Lebesgue measure on $H_N$. Fix positive real numbers $c_1 , \ldots , c_N$ and positive integers $\nu_1 , \ldots , \nu_N$. Then for $\RRe s > 0$ we have the following formula where $q = 1 / \sum^N_{i=1} \nu^{-1}_i$
\[
I := \int_{H_N} \Big( \sum^N_{i=1} c_i e^{-\nu_i x_i} \Big)^{\!-s} \! d\lambda = q (\nu_1 \cdots \nu_N)^{-1} \Big( \prod^N_{i=1} c^{q / \nu_i}_i \Big)^{\!-s} \! \Gamma (s)^{-1} \prod^N_{i=1} \Gamma (qs / \nu_i) \; .
\]
\end{uprop}

\begin{proof}
  First we show that the integral exists. Using lemma \ref{t431} and the fact that $\min (x_i) \le 0$ for $x \in H_N$, we find with $c = \min (c_i)$:
  \begin{equation}
    \label{eq:45}
    \sum^N_{i=1} c_i e^{-\nu_i x_i} \ge c \, e^{-\min (x_i)} \ge c \exp ((N-1)^{-1} \Betr{{x}}_{\infty}) \quad \mbox{for} \; x \in H_N \; .
  \end{equation}
Thus the function
\[
\Big( \sum^N_{i=1} c_i e^{-\nu_i x_i} \Big)^{-\RRe s} 
\]
is integrable over $H_N$. In order to evaluate the integral we recall the Mellin transform of a function $h$ on $\R^*_+$:
\[
(M h) (s) = \int^{\infty}_0 h (t) t^s \frac{dt}{t} \quad \mbox{for} \; \RRe s \ge 1
\]
and the convolution of two $L^1$-functions $h_1$ and $h_2$ on $\R^*_+$:
\[
(h_1 \ast h_2) (t) = \int^{\infty}_0 h_1 (t_1) h_2 (tt^{-1}_1) \frac{dt_1}{t_1} \; .
\]
For suitable $h_1$ and $h_2$ Fubini's theorem implies the basic formula
\[
M (h_1 \ast h_2) = (Mh_1) \cdot (Mh_2) \quad \mbox{for} \; \RRe s \ge 1 \; .
\]
For $t > 0$ let $d\mu$ be the image of Lebesgue measure under the exponential isomorphism:
\[
\{ x \in \R^N \tei \sum x_i = \log t \} \silo \{ (t_1 , \cdots , t_N) \in (\R^*_+)^N \tei t_1 \cdots t_N = t \} \; .
\]
The $N$-fold convolution of $L^1$-functions $h_1 , \ldots , h_N$ on $\R^*_+$ is given by the formula
\[
(h_1 \ast \ldots \ast h_N) (t) = \int_{t_1 \cdots t_N = t} h_1 (t_1) \cdots h_N (t_N) \, d\mu \; .
\]
Note that convolution is associative.

We may rewrite $I$ as follows
\[
I = \int_{t_1 \cdots t_N = 1} \Big( \sum^N_{i=1} c_i t^{\nu_i}_i \Big)^{-s} \, d\mu \; .
\]
Thus
\begin{eqnarray}
  \label{eq:46}
  \Gamma (s) \cdot I & = & \int^{\infty}_0 \Big( \int_{t_1 \cdots t_N = 1} \exp \Big( -t \sum^N_{i=1} c_i t^{\nu_i}_i \Big) d\mu \Big) t^s \frac{dt}{t}\\
& = & \int^{\infty}_0 \Big( \int_{t_1 \cdots t_N = t^{1/q}} \exp \Big( - \sum^N_{i=1} c_i t^{\nu_i}_i \Big) d\mu \Big) t^s \frac{dt}{t} \nonumber \\
& = & q M (e^{-c_1 t^{\nu_1}} \ast \ldots \ast e^{-c_N t^{\nu_N}}) (qs) \nonumber \\
& = & q M (e^{-c_1 t^{\nu_1}}) (qs) \cdots M (e^{-c_N t^{\nu_N}}) (qs) \nonumber \\
& = & q \prod^N_{i=1} \nu^{-1}_i c^{-qs /\nu_i}_i \Gamma (qs / \nu_i) \; . \nonumber
\end{eqnarray}
\end{proof}

We may now use the complex Stirling asymptotics
\begin{equation}
  \label{eq:47}
  \Gamma (s) \sim \sqrt{2\pi} e^{-s} e^{\left( s -\halb \right) \log s} \quad \mbox{for} \; \betr{s} \to \infty \; \mbox{in} \; - \pi < \arg s < \pi
\end{equation}
to draw the following consequence of proposition \ref{t432}.

\begin{ucor}
  \label{t433}
Let $k / \Q$ be a number field of degree $n$ with unit rank $r = r_1 + r_2 - 1 \ge 1$. Then we have the following asymptotic formula for $\RRe s \to \infty$, the integral being defined for $\RRe s > 0$:
\[
\int_{\sum_{v \tei \infty} x_v = 0} \Big( \sum_{v \, \real} e^{-2x_v} + 2 \sum_{v\,\complex} e^{-x_v} \Big)^{-s} d \lambda \sim \alpha_k s^{-r/2} n^{-s} \; .
\]
\end{ucor}

\begin{proof}
  Applying proposition \ref{t432} with $N = r_1 + r_2$ and the obvious choices of $c_i$'s and $\nu_i$'s the integral is seen to equal:
\[
n^{-1} 2^{1 - r_1} 2^{-2s r_2 /n} \Gamma (s)^{-1} \Gamma (s/n)^{r_1} \Gamma (2s / n)^{r_2} \; .
\]
Applying the Stirling asymptotics gives the result after some calculation.
\end{proof}

\begin{ucor}
  \label{t434}
Assumptions as in corollary \ref{t433}. For any $\varepsilon > 0$ set
\[
V_{\varepsilon} = \Big\{ x \in \bigoplus_{v\tei \infty} \R \tei \sum_{v \tei \infty} x_v = 0 \quad \mbox{and} \quad \Betr{{x}}_{\infty} < \varepsilon \Big\} \; .
\]
Then we have the asymptotic formula for $\RRe s \to \infty$:
\[
\int_{V_{\varepsilon}} \Big( \sum_{v\,\real} e^{-2x_v} + 2 \sum_{v\,\complex} e^{-x_v} \Big)^{-s} d\lambda \sim \alpha_k s^{-r/2} n^{-s} \; .
\]
\end{ucor}

\begin{proof}
  Set $f (x) = \sum_{v\,\real} e^{-2x_v} + 2\sum_{v\,\complex} e^{-x_v}$. For $x \in \bigoplus_{v\tei \infty} \R$ with $\sum_{v\tei \infty} x_v = 0$ we have by lemma \ref{t431} that:
  \begin{equation}
    \label{eq:48}
    f (x) \ge \exp (r^{-1} \Betr{{x}}_{\infty}) \; .
  \end{equation}
Choose $R \ge 2 r \log 2n$. For $\Betr{{x}}_{\infty} \ge R$ and $\alpha \ge 0$ we find
\[
\exp (- \alpha r^{-1} \Betr{{x}}_{\infty}) \le (2n)^{-\alpha} \exp \left( - \frac{\alpha}{2r} \Betr{{x}}_{\infty} \right) \; .
\]
For $\RRe s \ge 1$ this implies that
\begin{equation}
  \label{eq:49}
\Big| \int_{\sum x_v = 0 \atop \Betr{{x}}_{\infty} > R} f (x)^{-s} d\lambda \Big| \le \gamma (2n)^{-\RRe s}
\end{equation}
where
\[
\gamma = \int_{\sum x_v = 0 \atop \Betr{{x}}_{\infty} > R} \exp \left( - \frac{1}{2r} \Betr{{x}}_{\infty} \right) d\lambda < \infty \; .
\]
By the arithmetic-geometric mean inequality we see that in $\{ \sum x_v = 0 \}$ the function $f (x)$ has global minimum equal to $n$. We have $f (0) = n$ and $f (x) > n$ for all $x \neq 0$, c.f. the proof of theorem \ref{t42}, {\bf 1}. Choose $R \ge 2r \log 2n$ such that $R \ge \varepsilon$. Let $a_{\varepsilon , R}$ be the minimum of $f$ in the compact set $S_{\varepsilon , R}$ of $x$ with $\sum x_v = 0$ and $\varepsilon \le \Betr{{x}}_{\infty} \le R$. Then we have $a_{\varepsilon , R} > n$ and
\begin{equation}
  \label{eq:50}
  \Big| \int_{S_{\varepsilon,R}} f (x)^{-s} d\lambda \Big| \le \vol (S_{\varepsilon ,R}) a^{-\RRe s}_{\varepsilon , R} \quad \mbox{for} \; \RRe s \ge 0 \; .
\end{equation}
Using corollary \ref{t433} and the estimates (\ref{eq:49}) and (\ref{eq:50}) we find successively:
\[
\alpha_k s^{-r/2} n^{-s} \sim \int_{\sum x_v = 0} f (x)^{-s} d\lambda \sim \int_{\sum x_v = 0 \atop \Betr{x}_{\infty} \le R} f(x)^{-s} d\lambda \sim \int_{\sum x_v = 0 \atop \Betr{x}_{\infty} \le \varepsilon} f (x)^{-s} d \lambda \; .
\]
\end{proof}

We can now conclude the proof of theorem \ref{t43}. Let $\varepsilon > 0$ be so small that the image of $V_{\varepsilon}$ in $CH^1 (X_k)^0$ under the map $x \mapsto \Dh_x = [D_x]$ with $D_x = \sum_{v\tei \infty} x_v \cdot v$ is a homeomorphism onto its image $U_{\varepsilon}$. Moreover $\varepsilon > 0$ should be so small that $U_{\varepsilon}$ is contained in a neighborhood $U$ as in theorem \ref{t42}, {\bf 6}. Then we have
\begin{eqnarray}
  \label{eq:51}
  \int_{U_{\varepsilon}} \nu (\Dh) a (\Dh)^{-s} d^0 \Dh & = & \betr{\mu (k)} \int_{V_{\varepsilon}} \Big( \sum_{v\,\real} e^{-2x_v} + 2 \sum_{v\,\complex} e^{-x_v} \Big)^{-s} \, d\lambda \\
& \sim & \betr{\mu (k)} \alpha_k s^{-r/2} n^{-s} \quad \mbox{for} \; \RRe s \to \infty \nonumber
\end{eqnarray}
by corollary \ref{t434}. By theorem \ref{t42}, {\bf 1} and {\bf 2} (or {\bf 5}) the minimum $a_{U_{\varepsilon}}$ of $a = a (\Dh)$ on the compact set $CH^1 (X_k)^0 \ohne U_{\varepsilon}$ satisfies $a_{U_{\varepsilon}} > n$. Moreover $\nu = \nu (\Dh)$ is bounded, $\le d$ say. Together with the estimate
\[
\Big| \int_{CH^1 (X_k)^0 \ohne U_{\varepsilon}} \nu (\Dh) a (\Dh)^{-s} d^0 \Dh \Big| \le d \; \vol (CH^1 (X_k)^0) a^{-\RRe s}_{U_{\varepsilon}} \quad \mbox{for} \; \RRe s \ge 0
\]
the asymptotics (\ref{eq:51}) now imply the assertion of theorem \ref{t43}.
\end{proof}

\begin{remark}
  \label{t44}
\rm Using the asymptotic development of the $\Gamma$-function instead of (\ref{eq:47}) one can improve the assertion of theorem \ref{t43}. For example, the same proof shows that for any $\varphi \in (0 , \pi / 2)$ we have
\[
C (s) = \betr{\mu (k)} \alpha_k s^{-r/2} n^{-s} (1 + O (s^{-1})) \quad \mbox{as} \; \RRe s \to \infty
\]
in the angular domain $\betr{\arg s} < \varphi$. The $O$-constant depends on $\varphi$.
\end{remark}

\section{The two-variable zeta function as a regularized product}
\label{sec:5}

In this section we first review a theorem of Illies about the zeta-regularizability of entire functions of finite order.

We then apply his criterion to prove that $L (H^1 (X_k) , s,w)$ and $\zeta (X_k , s,w)$ are zeta-regularized as functions of $s$. 

There are many instances where one would like to give a sense to a non-convergent product of distinct non-zero complex numbers $a_{\nu}$ given with multiplicities $m_{\nu} \in \Z$. Sometimes the process of zeta regularization helps. Fix arguments $-\pi < \arg a_{\nu} \le \pi$ and assume that the Dirichlet series $D (u) = \sum m_{\nu} a^{-u}_{\nu}$ converges for $\RRe u \gg 0$ with a meromorphic continuation to $\RRe u > - \varepsilon$ for some $\varepsilon > 0$. If $D$ is holomorphic at $u = 0$ we may form the zeta-regularized product
\[
\rprod{(m_\nu)}{}\; a_{\nu} := \exp (-D' (0)) \; .
\]
If all $m_{\nu} = 1$, one sets $\rprod{}{} a_{\nu} = \rprod{(1)}{} a_{\nu}$.
In this way one obtains for example $\rprod{\infty}{\nu = 1} \nu = \sqrt{2\pi}$. For a finite sequence of $a_{\nu} , m_{\nu}$ the zeta-regularized product $\rprod{(m_{\nu})}{} a_{\nu}$ exists and equals the ordinary product $\prod a^{m_{\nu}}_{\nu}$.

For complex $s$ with $s \neq a_{\nu}$ for all $\nu$ one may ask whether $\rprod{(m_{\nu})}{} \; (s-a_{\nu})$ exists. In favourable instances it will define a meromorphic function in $\C$ whose zeroes and poles are precisely the numbers $a_{\nu}$ with their multiplicity $m_{\nu}$. On the other hand if we are given a meromorphic function $f (s)$ whose zeroes and poles are the numbers $a_{\nu}$ with multiplicity $m_{\nu}$ we may ask whether $\rprod{(m_{\nu})}{} \; (s - a_{\nu})$ exists and defines a meromorphic function in $\C$ and how it compares to $f (s)$. Sometimes it is also useful to introduce a scaling factor $\alpha > 0$ and compare $f (s)$ with $\rprod{(m_{\nu})}{} \; \alpha (s-a_{\nu})$. In the case where we have 
\[
f (s) = \rprod{(m_{\nu})}{} \; \alpha (s-a_{\nu}) \; ,
\]
the function $f$ is called ``$\alpha$-zeta regularized''. 

A much more thorough discussion of these problems and other regularization procedures ($\delta$-regularization) may be found in Illies' papers \cite{I1}, \cite{I2} and his references.

We now describe the precise technical result from Illies' work that we will use.

For $\varphi_1 , \varphi_2$ in $(0, \pi)$ define the open sets
\[
\Wh r_{\varphi_1 , \varphi_2} = \{ s \in \C^* \tei - \varphi_2 < \arg s < \varphi_1 \}
\]
and
\[
\Wh l_{\varphi_1 , \varphi_2} = \C^* \ohne \overline{\Wh r_{\varphi_1 , \varphi_2}} \; .
\]
A meromorphic function in $\C$ is said to be of finite order if it is the quotient of two entire functions of finite order.

\begin{theorem}[Illies]
  \label{t51}
Let $f$ be a meromorphic function of finite order in $\C$ such that almost all zeroes and poles lie in some $\Wh l_{\varphi_1 , \varphi_2}$. We assume that for some $0 < p \le \infty$ and any $p' < p$ we have 
\[
f (s) -1 = O (\betr{s}^{-p'}) \quad \mbox{in} \; \Wh r_{\varphi_1 , \varphi_2} \; \mbox{as} \; \betr{s} \to \infty \; .
\]
Then the following two assertions hold:\\
{\bf A} Setting $m (\rho) = \ord_{s = \rho} f (s)$, for any scaling factor $\alpha > 0$ the Dirichlet series 
\[
\xi (u,s) = \sum_{\rho \in \Wh l_{\varphi_1 , \varphi_2}} m (\rho) [\alpha (s- \rho)]^{-u}
\]
is uniformly convergent to a holomorphic function in $\RRe u \gg 0$ and $\betr{s} \ll 1$. Here we have chosen $- \pi < \arg (s - \rho) < \pi$ which is possible for small enough $\betr{s}$. The function $\xi (u,s)$ has a holomorphic continuation to any region of the form
\[
\{ \RRe u > -p \} \times G 
\]
where $G$ is an arbitary simply connected domain which does not contain zeroes or poles of $f$.\\
{\bf B} We have an equality of meromorphic functions in $\C$
\begin{eqnarray*}
  f (s) & = & \exp \Big( - \frac{\partial \xi}{\partial u} (0,s) \Big) \prod_{\rho \notin \Wh l_{\varphi_1, \varphi_2}} [\alpha (s-\rho)]^{m (\rho)} \\
& = & \rprod{(m (\rho))}{\rho} \; \alpha (s-\rho) \; .
\end{eqnarray*}
\end{theorem}

\begin{rem}
  According to {\bf B} the function $f$ equals the zeta-regularized determinant (scaled by $\alpha$) of its divisor. In fact $f$ is the $\delta$-regularized determinant of its divisor for any regularization sequence $\delta$ as in \cite{I2} Definition 3.4 but we do not need this stronger statement.
\end{rem}

\begin{proof}
  The result generalizes \cite{I2} Corollary 8.1 and is proved in the same way using \cite{I2} Theorem 5 and Proposition 3.3. The latter results are stated for the case where {\it all} zeroes and poles of $f$ lie in $\Wh l_{\varphi_1,\varphi_2}$. Using translation invariance of regularization as indicated in \cite{I2} Example 2) after Definition 4.1 gives the general case. In the thesis \cite{I1} more details can be found: Theorem 4.1 is a special case of \cite{I1} Korollar 2.7.1 and translation invariance is discussed in \cite{I1} Definition 2.2.3 and Korollar 2.2.4.
\end{proof}

We can now state our main theorem.

\begin{theorem}
  \label{t52}
For $k / \Q$ and any fixed complex number $w$ the functions \\
$\zeta (X_k , s,w)$ and $L (H^1 (X_k) , s,w)$ of $s$ are $\frac{1}{2\pi}$-zeta regularized.
\end{theorem}

In the function field case the corresponding but much simpler result is this

\begin{theorem}
  \label{t53}
Let $X / \F_q$ be a smooth projective and geometrically irreducible curve. Then for any fixed $w$ the entire function $P_X (q^{-s} , q^w)$ of $s$ and the meromorphic function $Z_X (q^{-s} , q^w)$ are $\alpha$-zeta regularized for any $\alpha > 0$.
\end{theorem}

\begin{rem}
  Comparing \ref{t52} and \ref{t53} we see that
\[
\zeta (X,s,w) := Z_X (q^{-s} , q^w) = (q^w-1)^{-1} q^{s (1-g)} \zeta_X (s,w)
\]
corresponds to
\[
\zeta (X_k , s,w) = w^{-1} \frac{2^{r_1/2}}{\betr{\mu (k)}} d^{-s/2}_k \zeta_{X_k} (s,w)
\]
in the following sense:  For every fixed $w$ both functions of $s$ are obtained by the process of $\frac{1}{2\pi}$-zeta regularization from the zeroes and poles of the analogous functions $(q^w - 1)^{-1} \zeta_X (s,w)$ and $w^{-1} \zeta_{X_k} (s,w)$. Note also that we have
\[
\zeta (X, s,1) = \zeta_X (s) \quad \mbox{and} \quad \zeta (X_k,s,1) = \hat{\zeta}_k (s) \; .
\]
\end{rem}

\begin{proofof}
  {\bf theorem \ref{t53}}
In view of formulas (\ref{eq:4}) and (\ref{eq:6}) this can be deduced from \cite{D2} 2.7 Lemma which evaluates $\rprod{}{\nu \in \Z} \alpha (s+\nu)$ for $\alpha \in \C^*$. Alternatively the theorem follows without difficulty from theorem \ref{t51}.
\end{proofof}

For the proof of theorem \ref{t52} we first need a refinement of the estimate given in Proposition \ref{t31} a).

\begin{lemma}
  \label{t54}
For any number field $k / \Q$ and every $R \ge 0$ there is a constant $c = c_k (R)$ such that setting 
\[
g (t, \Dh , w) = w^{-1} (k^0 (\Dh + \Dh_t)^w - 1- w \nu (\Dh) e^{-\pi t^{-2/n} a (\Dh)})
\]
we have
\[
\betr{g (t, \Dh , w)} \le c \exp (- \pi t^{-2/n} \min (2a (\Dh) , b (\Dh)))
\]
uniformly in $\Dh \in CH^1 (X_k)^0$ and $0 < t \le \sqrt{d_k}$ and $\betr{w} \le R$. 
\end{lemma}

\begin{proof}
  For $\Dh = [D]$ we have:
\[
k^0 (\Dh + \Dh_t) = \sum_{f \in I (D)} e^{-\pi t^{-2/n} \Betr{{f}}^2_D} \; .
\]
Hence
\[
\delta (t, \Dh) := k^0 (\Dh + \Dh_t) - 1 - \nu (\Dh) e^{-\pi t^{-2/n} a (\Dh)}
\]
is positive. We claim that there is a constant $\gamma$ depending only on $k$ such that for all $\Dh \in CH^1 (X_k)^0$ and $0 < t \le d^{1/2}_k$ we have the estimate
\begin{equation}
  \label{eq:52}
  0 < \delta (t , \Dh) \le \gamma e^{-\pi t^{-2/n} b (\Dh)} \; .
\end{equation}
This is seen as follows:
\begin{eqnarray*}
  \delta (t , \Dh) e^{\pi t^{-2/n} b (\Dh)} & = & \sum_{\Betr{{f}}^2_D \ge b (\Dh)} e^{-\pi t^{-2/n} (\Betr{{f}}^2_D - b (\Dh))}\\
& \le & \sum_{\Betr{{f}}^2_D \ge b (\Dh)} e^{-\pi d^{-1/n}_k (\Betr{{f}}^2_D - b (\Dh))} \quad \mbox{since} \; t \le \sqrt{d_k} \\
& \le & e^{\pi d^{-1/n}_k b (\Dh)} k^0 (\Dh + \Dh_{\sqrt{d_k}}) = f (\Dh) \; .
\end{eqnarray*}
Hence for $\gamma$ we may choose the supremum of the bounded function $f$ on $CH^1 (X_k)^0$, c.f. Proposition \ref{t41}.

Since the left hand side of the estimate in lemma \ref{t54} is bounded and since $a = a (\Dh)$ is bounded on $CH^1 (X_k)^0$ it sufficed to prove the desired estimate for all $0 < t \le \varepsilon$, where $\varepsilon > 0$ is small. We choose $0 < \varepsilon \le \sqrt{d_k}$ such that for all $0 < t \le \varepsilon$ and $\Dh \in CH^1 (X_k)^0$ we have
\[
0 < x = x (t, \Dh) = k^0 (\Dh + \Dh_t) - 1 \le 1/2 \; .
\]
This is possible by (\ref{eq:24}) or (\ref{eq:52}). We may assume that $R \ge 1$. Using inequality (\ref{eq:26}), we find:
\begin{eqnarray*}
  k^0 (\Dh + \Dh_t)^w & = & (1 + x)^w = 1 + wx + w x^2 \vartheta\\
& = & 1 + w \nu (\Dh) e^{-\pi t^{-2/n} a (\Dh)} + w \psi
\end{eqnarray*}
where
\[
\betr{\vartheta} = \betr{\vartheta (t , w, \Dh)} \le e^{2R}
\]
and
\[
\psi = \delta (t , \Dh) + x^2 \vartheta \; .
\]
From (\ref{eq:52}) we get
\begin{eqnarray*}
\betr{\psi} & \le & \gamma e^{-\pi t^{-2/n} b (\Dh)} + e^{2R} e^{-2\pi t^{-2/n} a (\Dh)} (\nu (\Dh) + \gamma e^{-\pi t^{-2/n} (b (\Dh) - a (\Dh))} )^2\\
& \le & \gamma e^{-\pi t^{-2/n} b (\Dh)} + e^{2R} e^{-2 \pi t^{-2/n} a (\Dh)} (\nu (\Dh) + \gamma)^2 \; .
\end{eqnarray*}
This gives the required estimate in $0 < t \le \varepsilon$ for $c = e^{2R} \max_{\Dh} (\nu (\Dh) + \gamma)$. 
\end{proof}

\begin{proofof}
  {\bf theorem \ref{t52}}
According to lemma \ref{t33} we may write the function $\zeta (X_k ,s,w)$ as follows in the region $\RRe s > \RRe w , \RRe s > 0$
\begin{eqnarray*}
  \zeta (X_k ,s,w) & = & w^{-1} \frac{2^{r_1/2}}{\betr{\mu (k)}} \int_{CH^1 (X_k)} (k^0 (\Dh)^w - 1) \Nh \Dh^{-s} d \Dh \\
& = & \frac{2^{r_1/2}}{\betr{\mu (k)}} \int^{\infty}_0 \int_{CH^1 (X_k)^0} w^{-1} (k^0 (\Dh + \Dh_t)^w - 1) d^0 \Dh \; t^{-s} \frac{dt}{t} \\
& = & \frac{2^{r_1/2}}{\betr{\mu (k)}}\int^{\infty}_0 \int_{CH^1 (X_k)^0} (\nu (\Dh) e^{-\pi t^{-2/n} a (\Dh)} + g (t , \Dh , w)) d^0 \Dh \; t^{-s} \frac{dt}{t} \; .
\end{eqnarray*}
The first term leads to the following meromorphic function in $\C$
\begin{eqnarray*}
A (s) & := & \int^{\infty}_0 \int_{CH^1 (X_k)^0} \nu (\Dh) e^{-\pi t^{-2/n} a (\Dh)} d^0 \Dh \; t^{-s} \frac{dt}{t} \\
 & = & \frac{n}{2} \pi^{-\frac{ns}{2}} \Gamma \Big( \frac{ns}{2} \Big) \int_{CH^1 (X_k)^0} \nu (\Dh) a (\Dh)^{-\frac{ns}{2}} d^0 \Dh \; .
\end{eqnarray*}
Setting
\[
f_w (s) = 1 + A (s)^{-1} \int^{\infty}_0 \int_{CH^1 (X_k)^0} g (t , \Dh , w) d^0 \Dh \; t^{-s} \frac{dt}{t}
\]
we may write
\begin{equation}
  \label{eq:53}
  \zeta (X_k ,s,w) = \frac{2^{r_1/2}}{\betr{\mu (k)}} A (s) f_w (s) \; .
\end{equation}
We will now show that for any fixed $w \in \C$ and $\alpha > 0$ the function $f_w (s)$ is $\alpha$-zeta regularized and that its $\xi (u,s)$-function in the sense of theorem \ref{t51} {\bf A} has a holomorphic continuation to any region $\C \times G$ where $G \subset \C$ is any simply connected domain disjoint from the zeroes and poles of $f_w$. First note that $f_w$ is meromorphic of finite order $(\le 1)$ since this is true for $s \mapsto \zeta (X_k ,s,w)$ by proposition 3.5 and clear for $\frac{2^{r_1/2}}{\betr{\mu (k)}} A (s)$.

Let $0 < \varphi < \pi/2$ be any angle such that 
\[
\varphi \; \tan \varphi < \log \min (2 , b_{\inf} / n) \; .
\]
Note here that because of theorem \ref{t42}, {\bf 1} and {\bf 7} we have $b_{\inf} > a_{\min} = n$. 

We will now show that for any $p' > 0$ and every $w \in \C$ we have the estimate
\begin{equation}
  \label{eq:54}
  f_w (s) - 1 = O (\betr{s}^{-p'}) \quad \mbox{in} \; \overline{\Wh r_{\varphi , \varphi}} \;\mbox{as} \; \betr{s} \to \infty \;  .
\end{equation}
It follows in particular that $f_w$ has only finitely many zeroes or poles in $\overline{\Wh r_{\varphi , \varphi}}$. Hence the conclusions of Illies' theorem 5.1 apply to $f_w$ and $\varphi_1 = \varphi_2 = \varphi$ with $p = \infty$. 

In order to prove (\ref{eq:54}) we have to show the following estimates for any fixed $w$ and $p' > 0$ as $\RRe s \to \infty$ in $\overline{\Wh r_{\varphi,\varphi}}$:
\begin{equation}
  \label{eq:55}
  A (s)^{-1} \int^{\sqrt{d}_k}_0 \int_{CH^1 (X_k)^0} g (t , \Dh , w) d^0 \Dh \; t^{-s} \frac{dt}{t} = O (\betr{s}^{-p'})
\end{equation}
and
\begin{equation}
  \label{eq:56}
  A (s)^{-1} \int^{\infty}_{\sqrt{d}_k} \int_{CH^1 (X_k)^0} g (t , \Dh ,w ) d^0 \Dh \; t^{-s} \frac{dt}{t} = O (\betr{s}^{-p'}) \; .
\end{equation}
We begin with (\ref{eq:55}). By lemma \ref{t54} we have for $0 < t \le \sqrt{d}_k$ and \\
$\Dh \in CH^1 (X_k)^0$:
\begin{eqnarray*}
  \betr{g (t, \Dh , w)} & \le & c \exp (-\pi t^{-2/n} \min (2a (\Dh) , b (\Dh))) \\
& \le & c \exp (- \pi t^{-2/n} \min (2a_{\min} , b_{\inf})) \; .
\end{eqnarray*}
Hence there are constants $c_1 , c_2$ depending on $w$ such that for $\RRe s \to \infty$ we have
\begin{eqnarray}
  \label{eq:57}
  \betr{A (s)^{-1} \int^{\sqrt{d}_k}_0 \cdots} & \le & c_1 \betr{A (s)}^{-1} (\pi \min (2a_{\min} , b_{\inf}))^{- \frac{n \RRe s}{2}} \Gamma \Big( \frac{n \RRe s}{2} \Big) \nonumber \\
& \le & c_2 \frac{\Gamma \left( \frac{n \RRe s}{2} \right) }{\betr{\Gamma \left( \frac{ns}{2} \right)}} \betr{s}^{r/2} \Big( \frac{n}{\min (2n , b_{\inf})} \Big)^{\frac{n \RRe s}{2}} \; .
\end{eqnarray}
For the second inequality we have used theorem \ref{t43} which was the main result of section \ref{sec:4} and the fact that $a_{\min} = n$. Now the Stirling asymptotics (\ref{eq:47}) shows that for any $\varphi \in (0, \pi / 2)$ there is a constant $c_{\varphi}$ such that for all $z \in \C$ with $\RRe z \ge 1/2$ and $z \in \overline{\Wh r_{\varphi , \varphi}}$ we have the estimate
\begin{equation}
  \label{eq:58}
  \frac{\Gamma (\RRe z)}{\betr{\Gamma (z)}} \le c_{\varphi} \exp ((\RRe z) \varphi \; \tan \varphi) \; .
\end{equation}
Namely, setting $z = r e^{i\alpha}$ with $\betr{\alpha} \le \varphi$ we have as $r \to \infty$
\begin{eqnarray*}
  \frac{\Gamma (\RRe z)}{\betr{\Gamma (z)}} & \sim & \exp ((r \cos \alpha - \halb) \log \cos \alpha + r \alpha \sin \alpha) \\
& \le & \exp (r \alpha \sin \alpha) = \exp ((\RRe z) \alpha \tan \alpha) \le \exp (( \RRe z) \varphi \tan \varphi) \; .
\end{eqnarray*}
Thus we can proceed with the estimate (\ref{eq:57}), obtaining
\begin{equation}
  \label{eq:59}
  \betr{A (s)^{-1} \int^{\sqrt{d}_k}_0 \cdots} \le c_3 \betr{s}^{r/2} \exp \Big( \frac{n \RRe s}{2} ( \varphi \tan \varphi - \log \min (2 , b_{\inf} / n)) \Big) \; .
\end{equation}
By our choice of $\varphi$, the second term converges exponentially fast to zero as $\RRe s \to \infty$. Since $\betr{s} \le (\RRe s) (1 + \tan^2 \varphi)^{1/2}$ in $\overline{\Wh r_{\varphi , \varphi}}$ the estimate (\ref{eq:55}) follows for all $p' > 0$.

Next, we prove (\ref{eq:56}). We have
\begin{eqnarray*}
\lefteqn{  \betr{g (t, \Dh,w)} \le \betr{w^{-1} (k^0 (\Dh + \Dh_t)^w-1)} + \nu (\Dh) e^{-\pi t^{-2/n} a (\Dh)}} \\
& \le & \betr{w^{-1} (k^0 (\Dh + \Dh_t)^w - t^w d^{-w/2}_k)} + \betr{w^{-1} (t^w d^{-w/2}_k - 1)} + \nu (\Dh) e^{- \pi t^{-2/n} a (\Dh)} \; .
\end{eqnarray*}
For $\Dh \in CH^1 (X_k)^0$ and $t \ge \sqrt{d}_k$ it follows from proposition \ref{t31} b) and the boundedness of $\nu$ on $CH^1 (X_k)^0$ that we have
\begin{eqnarray*}
  \betr{g (t, \Dh , w)} & \le & c_4 \exp (- \alpha t^{2/n}) + \betr{w^{-1} (t^w d^{-w/2}_k - 1)} + c_5 e^{-\pi n t^{-2/n}} \\
& \le & c_6 t^M \; .
\end{eqnarray*}
Here $M = \max (\RRe w , 1)$ will do, and the constants $c_i$ depend on $w$. Observe that for $w = 0$ the middle term becomes a logarithm in $t$ which is absorbed in $t^M$ since $M \ge 1 > 0$. 

Using the estimate and theorem \ref{t43} we find for $\RRe s > M$ in $\overline{\Wh r_{\varphi , \varphi}}$:
\begin{eqnarray*}
  \betr{A (s)^{-1} \int^{\infty}_{\sqrt{d}_k} \cdots } & \le & c_7 \betr{A (s)}^{-1} \betr{\int^{\infty}_{\sqrt{d}_k} t^{M-s} \frac{dt}{t}} \\
& \le & c_8 \betr{\Gamma \Big( \frac{ns}{2} \Big)}^{-1} \betr{s}^{r/2} (\pi n)^{\frac{n \RRe s}{2}} \betr{M-s}^{-1} d^{-\frac{\RRe s}{2}}_k \\
 & \le & c_9 \betr{e^{-\frac{ns}{2} \log \frac{ns}{2}}} \betr{s}^{(r+1)/2} (\pi en)^{\frac{n \RRe s}{2}} \betr{M-s}^{-1} d^{-\frac{\RRe s}{2}}_k
\end{eqnarray*}
by the Stirling asymptotics. Together with the estimate
\[
\betr{e^{-\frac{ns}{2} \log \frac{ns}{2}}} \le e^{- \frac{n\RRe s}{2} \log \betr{\frac{ns}{2}}} e^{\frac{n \RRe s}{2} \varphi \tan \varphi}
\]
this implies the desired estimate (\ref{eq:56}) for any $p' > 0$.
\medskip

Having thus proved (\ref{eq:54}), theorem \ref{t51} implies that $f_w$ is $\alpha$-zeta regularized for any $\alpha > 0$. Now, equation (\ref{eq:53}) together with formula (\ref{eq:39}) gives
\begin{equation}
  \label{eq:60}
  \zeta (X_k , s,w) = \hat{\zeta}_k (s) \frac{f_w (s)}{f_1 (s)} \; .
\end{equation}
It has been known for quite some time that $\hat{\zeta}_k (s)$ is $\frac{1}{2\pi}$-zeta regularized and there are different ways to see this. For example $\zeta_k (s)$ is $\alpha$-zeta regularized for any positive $\alpha > 0$ by \cite{I2} corollary 8.1 which applies to very general Dirichlet series. Furthermore the $\Gamma$-factors $\Gamma_{\R} (s)$ and $\Gamma_{\C} (s)$ are $\frac{1}{2\pi}$-zeta regularized as follows from a formula essentially due to Lerch, \cite{D2} (2.7.1):
\begin{equation}
\label{eq:61}
\rprod{\infty}{\nu=0} \alpha (z + \nu) = \alpha^{\halb - z} \Big( \frac{1}{\sqrt{2\pi}} \Gamma (z) \Big)^{-1} \; .
\end{equation}
It follows from (\ref{eq:60}) that $s \mapsto \zeta (X_k ,s,w)$ is $\frac{1}{2\pi}$-zeta regularized for any $w$. Hence theorem \ref{t52} is proved.
\end{proofof}

Incidentally, we may deduce the following corollary from the proof of \ref{t52}:

\begin{cor}
  \label{t55}
For any number field $k / \Q$ the entire function
\[
\tilde{C} (s) = 2^{r_1/2} \sqrt{n/2} \int_{CH^1 (X_k)^0} \frac{\nu (\Dh)}{\betr{\mu (k)}} \Big( \frac{a (\Dh)}{n} \Big)^{-\frac{ns}{2}} d^0 \Dh
\]
is $\frac{1}{2\pi}$-zeta regularized.
\end{cor}
\medskip

\noindent {\bf 1. Proof}  According to the above and formula (\ref{eq:53}) for $w = 1$, the function
  \begin{equation}
    \label{eq:62}
    \frac{\hat{\zeta}_k (s)}{f_1 (s)} = \frac{2^{r_1/2}}{\betr{\mu (k)}} \frac{n}{2} \pi^{-\frac{ns}{2}} \Gamma \Big( \frac{ns}{2} \Big) \int_{CH^1 (X_k)^0} \nu (\Dh) a (\Dh)^{- \frac{ns}{2}} d^0 \Dh
  \end{equation}
is $\frac{1}{2\pi}$-zeta regularized. It follows from formula (\ref{eq:61}) that we have
\[
\pi^{-\frac{ns}{2}} \Gamma \Big( \frac{ns}{2} \Big) = n^{\frac{ns}{2}} \sqrt{2/n} \Big( \rprod{\infty}{\nu=0} \frac{1}{2\pi} (s + \frac{2\nu}{n} ) \Big)^{-1} \; .
\]
Up to a $\frac{1}{2\pi}$-zeta regularized function the term $\pi^{-\frac{ns}{2}} \Gamma (ns/2)$ in formula (\ref{eq:62}) can therefore be replaced by $n^{\frac{ns}{2}} \sqrt{2/n}$. This gives the assertion. 
\beweisende
\bigskip

\noindent {\bf 2. Proof} It follows from remark \ref{t44} that for any $\varphi \in (0 , \pi/2)$ we have
\[
(s / 2 \pi)^{r / 2} \tilde{C} (s) = 1 + O (s^{-1}) \quad \mbox{in} \; \overline{\Wh r_{\varphi , \varphi}} \; \mbox{as} \; \RRe s \to \infty \; .
\]
By theorem \ref{t51} this function is therefore even $\alpha$-zeta regularized for every $\alpha > 0$. \beweisende
\bigskip

Let us check the corollary for $k = \Q$ and $k$ imaginary quadratic. For $k = \Q$ the function equals $1$ which is regularized. For $k$ imaginary quadratic the function reduces to the integral, which in this case is a finite Dirichlet series over ideal classes. Because of $\nu (0) = \betr{\mu (k)}$ and $a (0) = n$ this Dirichlet series starts with a constant term $1$. Now \cite{I2} Corollary 8.1, resp. its proof shows that such a finite Dirichlet series is $\alpha$-zeta regularized for any $\alpha > 0$.

\section{The cohomological motivation}
\label{sec:6}

In this section we explain how theorem \ref{t52} fits into the speculative cohomological setting of \cite{D3}. For every number field $k / \Q$ there should exist complex topological cohomology spaces $H^i (X_k , \Ch)$ together with an $\R$-action $\Phi^t$. The infinitesimal generator $\Theta$ of this $\R$-action should exist. We expect that
\[
H^0 (X_k , \Ch) = \C \quad \mbox{with} \; \Theta = 0 
\]
and
\[
H^2 (X_k , \Ch) \silo \C \quad \mbox{with} \; \Theta = \id \; .
\]
The space $H^1 (X_k , \Ch)$ should be infinite dimensional and decompose in a suitable sense into the eigenspaces of $\Theta$, the eigenvalues being the zeroes of $\hat{\zeta}_k (s)$. In degrees greater than two the cohomologies should vanish. 

The zeta-regularized determinant $\det_{\infty} (\varphi)$ of a diagonalizable operator $\varphi$ is defined as the zeta-regularized product of its eigenvalues with their multiplicities. See \cite{D2} for more precise definitions. The relation between $\hat{\zeta}_k (s)$ and cohomology is expected to be:
\begin{equation}
  \label{eq:63}
  \hat{\zeta}_k (s) = \prod^2_{i=0} \ddet_{\infty} \Big( \frac{1}{2\pi} (s \cdot \id - \Theta) \tei H^i (X_k , \Ch) \Big)^{(-1)^{i+1}} \; .
\end{equation}
From this and the above it follows that we would have
\begin{eqnarray}
  \label{eq:64}
  L (H^1 (X_k) , s) & := & \frac{s}{2\pi} \frac{s-1}{2\pi} \hat{\zeta}_k (s) \nonumber \\
 & = & \ddet_{\infty} \Big( \frac{1}{2\pi} (s \cdot \id - \Theta) \tei H^1 (X_k, \Ch) \Big) \; .
\end{eqnarray}
Formulas (\ref{eq:63}) and (\ref{eq:64}) would imply in particular that $\hat{\zeta}_k (s)$ and $L (H^1 (X_k) , s)$ are $\frac{1}{2\pi}$-zeta regularized and this turned out to be true, \cite{D1} \S\,4, \cite{SchS}, \cite{JL}, \cite{I1}.

How to incorporate the two-variable zeta function into this picture? One natural idea is to assume that there is an operator $\Theta_w$ on $H^{\hullet} (X_k , \Ch)$ for every $w \in \C$ deforming $\Theta_1 = \Theta$ and such that the two variable zeta-function equals
\begin{equation}
  \label{eq:65}
  \prod^2_{i=0} \ddet_{\infty} \Big( \frac{1}{2\pi} (s \cdot \id - \Theta_w) \tei H^i (X_k , \Ch)\Big)^{(-1)^{i+1}} \; .
\end{equation}
The equation (\ref{eq:32})
\[
\hat{\zeta}_k (s) = \frac{2^{r_1/2}}{\betr{\mu (k)}} d^{-s/2}_k \zeta_{X_k} (s,1)
\]
and formula (\ref{eq:63}) suggest that the function
\begin{equation}
  \label{eq:66}
  \frac{2^{r_1/2}}{\betr{\mu (k)}} d^{-s/2}_k \zeta_{X_k} (s,w)
\end{equation}
might be equal to (\ref{eq:65}). However the function (\ref{eq:66}) is identically zero for $w = 0$ and this is incompatible with (\ref{eq:65}). Namely the zeroes of (\ref{eq:65}) come from factors of the form $\frac{1}{2\pi} (s-\lambda)$ where $\lambda \in \spec (\Theta_w)$ if the zeta-regularized products exist in the sense recalled in section \ref{sec:5}. The easiest modification of (\ref{eq:66}) which takes this point into account is to consider instead of (\ref{eq:66}) the function
\[
w^{-1} \frac{2^{r_1/2}}{\betr{\mu (k)}} d^{-s/2}_k \zeta_{X_k} (s,w)
\]
i.e. $\zeta (X_k , s,w)$. Thus the following equation is suggested
\begin{equation}
  \label{eq:67}
  \zeta (X_k , s,w) = \prod^2_{i=0} \ddet_{\infty} \Big( \frac{1}{2\pi} (s \cdot \id - \Theta_w) \tei H^i (X_k , \Ch) \Big)^{(-1)^{i+1}} \; .
\end{equation}
It would imply that $\zeta (X_k , s,w)$ is $\frac{1}{2\pi}$-zeta regularized. This was proved in theorem \ref{t52}.

The poles of $s \mapsto \zeta (X_k ,s,w)$ lie at $s = 0$ and $s = w$. For $w \neq 0$ they have order one. For $w = 0$ there is a double pole at $s = 0$. According to (\ref{eq:67}) the poles of $\zeta (X_k ,s,w)$ are accounted for by the eigenvalues of $\Theta_w$ on $H^0 (X_k, \Ch)$ and $H^2 (X_k , \Ch)$. On $H^0 (X_k , \Ch) = \C$ it is natural to expect $\Theta_w = 0$ for all $w$. It follows that on $H^2 (X_k , \Ch) \cong \C$ we {\it must} have $\Theta_w = w \cdot \id$. Then (\ref{eq:67}) implies the formula
\[
\frac{s}{2\pi} \frac{s-w}{2\pi} \zeta (X_k , s,w) = \ddet_{\infty} \Big( \frac{1}{2\pi} (s \cdot \id - \Theta_w) \tei H^1 (X_k , \Ch) \Big) \; .
\]
This is the reason why we denoted the left hand side by $L (H^1 (X_k) , s , w)$ in definition \ref{t34}.

Having explained the motivation behind theorem \ref{t52} let us discuss the speculative formula (\ref{eq:67}) a little further. The functional equation (\ref{eq:39}) for $\zeta (X_k ,s,w)$ says in particular that $\rho \mapsto w-\rho$ is an involution on the set of zeroes resp. poles of $s \mapsto \zeta (X_k ,s,w)$. Under (\ref{eq:67}) this is compatible with the expected Poincar\'e duality
\[
\cup : H^i (X_k , \Ch) \times H^{2-i} (X_k , \Ch) \longrightarrow H^2 (X_k , \Ch) \cong \C
\]
if we assume that $\Theta_w$ is a derivation with respect to $\cup$-product. It looks like $\Theta_w$ was the infinitesimal generator of an $\R$-action $\Phi^t_w$ on cohomology which respects cup product. It could be interesting to check whether there is a symplectic structure in the distribution of the low lying zeroes of $s \mapsto L (H^1 (k) , s,w)$ as in the work of Katz and Sarnak \cite{S}.

In contrast to $\Theta$ the operators $\Theta_w$ for $w < 0$ will not commute with the Hodge $*$-operator as in \cite{D3} \S\,3 since this would force the zeroes of $\zeta (X_k , s,w)$ to lie on the line $\RRe s = \frac{w}{2}$ which is not the case for $w < 0$ by the investigations of Lagarias and Rains, \cite{LR} \S\,7.

From calculations in the function field case, I do not expect the operators $\Theta_w$ for different $w$ to commute. One possibility seems to be that $[\Theta_{w_1} , \Theta_{w_2}] = (w_1 - w_2) \id$. 

\newpage
\noindent
Mathematisches Institut\\
Westf. Wilhelms-Universit\"at\\
Einsteinstr. 62\\
48149 M\"unster\\
Germany\\
deninge@math.uni-muenster.de
%\vspace*{0.5cm}

%\hspace*{\fill} {\tiny Stand: \today}\\
%\input{bild}

\begin{thebibliography}{9999}
\bibitem[D1]{D1}C. Deninger, Local $L$-factors of motives and regularized determinants. Invent. Math. {\bf 107} (1992), 135--150
\bibitem[D2]{D2}C. Deninger, Motivic $L$-functions and regularized determinants. Proc. Symp. Pure Math. {\bf 55}, 1 (1994), 707--743
\bibitem[D3]{D3}C. Deninger, Some analogies between number theory and dynamical systems on foliated spaces. Doc. Math. J. DMV. Extra Volume ICM I (1998), 23--46
%\bibitem[D4]{D4}C. Deninger, Number theory and dynamical systems on foliated spaces. Jber. d. Dt. Math. Verein. {\bf 103} (2001), 79--100
\bibitem[F]{F}P. Francini, The size function $h^0$ for quadratic number fields. Journal de Th\'eorie des Nombres de Bordaux {\bf 13} (2001), 125--135
\bibitem[GS]{GS}G. van der Geer and R. Schoof, Effectivity of Arakelov divisors and the theta divisor of a number field. Selecta Math., New Series {\bf 6} (2000), 377--398
\bibitem[I1]{I1}G. Illies, Regularisierte Produkte, Spurformeln und Cram\'er-Funktionen. Schriftenreihe des Math. Inst. Univ. M\"unster, 3. series, Heft {\bf 22} (1998)
\bibitem[I2]{I2}G. Illies, Regularized products and determinants. Commun. Math. Phys. {\bf 220} (2001), 69--94
\bibitem[JL]{JL}J. Jorgenson, S. Lang, Basic analysis of regularized series and products, LNM {\bf 1564}, Springer 1994
%\bibitem[K]{K}K. Knopp, Theorie und Anwendung der unendlichen Reihen. Grundlehren Band 2, Springer 1964
\bibitem[LR]{LR}J.C. Lagarias, E. Rains, On a two-variable zeta function for number fields. ArXiv: math.NT/0104176v5, 7 July 2002
\bibitem[N]{N}N. Naumann, On the irreducibility of the two variable zeta-function for curves over finite fields. ArXiv:math.AG/0209092, 2002
\bibitem[P]{P}R. Pellikaan, On special divisors and the two variable zeta function of algebraic curves over finite fields. In: R. Pellikaan, M. Perret, S.G. Vladut, eds.: Arithmetic, Geometry and Coding theory, Walter de Gruyter, Berlin 1996, 175--184
\bibitem[S]{S}P. Sarnak, $L$-functions. Doc. Math. J. DMV. Extra volume ICM I (1998), 453--465
\bibitem[SchS]{SchS}M. Schr\"oter, C. Soul\'e, On a result of Deninger concerning Riemann's zeta function. Proc. Symp. Pure Math. {\bf 55}, 1 (1994), 745--747
\end{thebibliography}
\end{document}